\documentstyle[12pt,amssymb]{article}

\begin{document}

\begin{center}
{\bf INTEGRALITY OF NEARLY (HOLOMORPHIC) SIEGEL MODULAR FORMS} 
\end{center}

\begin{center}
By {\sc Takashi Ichikawa} 
\end{center}

\centerline{\rule[1mm]{30mm}{0.05mm}}
\vspace{2ex}

\noindent
\begin{small}
{\it Abstract.} 
In order to considering the integrality of nearly holomorphic (vector-valued) 
Siegel modular forms, 
we introduce nearly Siegel modular forms and study their integrality. 
We show that the integrality of nearly Siegel modular forms 
in terms of their Fourier expansion implies the integrality of their CM values. 
Furthermore, we show that there exists a one-to-one correspondence 
between integral nearly Siegel modular forms and integral nearly holomorphic ones. 
By these results,  the integrality of  CM values holds for 
nearly holomorphic Siegel modular forms 
and for nearly overconvergent $p$-adic Siegel modular forms. 
\end{small}
\vspace{3ex}

{\bf 1. Introduction.} 
Shimura's work (summarized in \cite{Sh}) on nearly holomorphic 
vector-valued modular forms of several variables 
is a fundamental theory 
to study the arithmetic of modular forms and $L$-functions. 
For example, he applied his theory to show the algebraicity of 
special values of the $L$-functions associated with Hecke characters and 
modular forms. 
In \cite{Kt1, Kt2, Kt3}, 
Katz gave a perspective of constructing $p$-adic version of Shimura's theory, 
and realized this perspective in the elliptic and Hilbert modular cases 
based on the rationality of CM values for arithmetic modular forms. 
Recently, similar results are given by Eischen \cite{E1,E2} for the unitary modular case 
(see also \cite{HLS}). 

The purpose of this paper is to study the integrality and $p$-adic property of 
CM values for nearly holomorphic (vector-valued) Siegel modular forms 
(of degree $>1$ and level $\geq 3$). 
For this purpose, 
it is natural to introduce {\it nearly Siegel modular forms} which are 
an algebraic version of nearly holomorphic Siegel modular forms. 
Nearly Siegel modular forms are considered in Darmon-Rotger \cite[Section 2]{DR} 
and Urban \cite{U} in the elliptic modular case, 
and are defined as global sections a vector bundle arising from the de Rham bundle 
over the Shimura model of a Siegel modular variety. 
The advantage point to considering nearly Siegel modular forms is that 
we can define such modular forms over a ring and study their integrality. 
We construct their arithmetic Fourier expansion satisfying the $q$-expansion principle. 
Therefore, 
one can determine the integrality of nearly Siegel modular forms 
by their Fourier expansion.  

Furthermore, 
we consider the analytic and $p$-adic realizations of nearly Siegel modular forms. 
The Hodge decomposition of the de Rham bundle gives the analytic realization of 
nearly Siegel modular forms as nearly holomorphic Siegel modular forms, 
and the unit root splitting (cf. \cite[1.11]{Kt3}) gives the $p$-adic realization of 
nearly Siegel modular forms as $p$-adic Siegel modular forms. 
The key result of this paper states, roughly speaking, the following: 
\vspace{2ex}

{\sc Theorem} (see Theorem 3.4 for the precise statement). 
\begin{it} 
The analytic realization map gives an isomorphism between the spaces 
of integral nearly Siegel modular forms and 
of integral nearly holomorphic Siegel modular forms with same weight. 
\end{it} 
\vspace{2ex}

By this theorem and the theory of nearly Siegel modular forms, 
one can see that the $q$-expansion principle and the integrality of CM values hold 
for nearly holomorphic Siegel modular forms. 
Further, we show that the $p$-adic realization map gives an injection 
from the space of nearly holomorphic Siegel modular forms into 
that of {\it nearly overconvergent} $p$-adic Siegel modular forms. 
It seems interesting to study the classicity problem, 
namely characterizing the image of the $p$-adic realization map 
based on results of Andreatta-Iovita-Pilloni \cite{AIP}. 
\vspace{2ex}

The organization of this paper is as follows. 

In Section 2, 
we define nearly Siegel modular forms as global sections of 
an automorphic de Rham bundle which is an extension of an automorphic Hodge bundle  on the Shimura model of a Siegel modular variety. 
We show fundamental properties of the space of nearly Siegel modular forms 
of fixed weight, 
for example that this space is finitely generated and has 
the arithmetic Fourier expansion satisfying the $q$-expansion principle. 

In Section 3, 
we give the analytic realization of nearly Siegel modular forms as 
nearly holomorphic Siegel modular forms. 
According to Shimura's theory, 
we define the integrality of nearly holomorphic Siegel modular forms 
by their Fourier expansion, 
and show the above key result with application to their integrality of CM values. 

In Section 4, 
we give the $p$-adic realization of nearly Siegel modular forms as 
$p$-adic Siegel modular forms. 
We show that this realization map becomes the composite of 
the analytic realization map and a linear map preserving $p$-ordinary CM values 
between the spaces of nearly holomorphic and overconvergent Siegel modular forms. 
Furthermore, 
we give examples of Siegel-Eisenstein series. 
\vspace{3ex}

{\bf 2. Nearly modular forms} 
\vspace{2ex}

{\bf 2.1. Representation of classical groups.} \ 
Let $V$ be a $2g$-dimensional vector space with symplectic form, 
and $W$ be its anisotropic subspace of dimension $g$. 
Then $GL_{g} = GL(W)$ is a general linear group of rank $g$ which is contained in 
a symplectic group $Sp_{2g} = Sp(V)$ of rank $g$ as 
$$
GL_{g} \cong 
\left\{ \left. \left( \begin{array}{cc} A & O \\ O & \mbox{}^{t}A^{-1} \end{array} \right) 
\in Sp_{2g} \ \right| \ A \in GL_{g} \right\}. 
$$
Let $B_{g}$ be the Borel subgroup of $GL_{g}$ consisting of upper-triangular matrices, 
and $B_{2g}$ denote the Borel subgroup of $Sp_{2g}$ given by 
$$
\left\{ \left. \left( \begin{array}{cc} A & * \\ O & \mbox{}^{t}A^{-1} \end{array} \right) 
\in Sp_{2g} \ \right| \ A \in B_{g} \right\}. 
$$
Then the maximal torus $T_{g} \subset B_{g}$ of $GL_{g}$ becomes that of $Sp_{2g}$, 
and ${\mathbb Z}^{g}$ is identified with the group $X(T_{g})$ of characters of $T_{g}$ as 
$$
\left( \begin{array}{ccc}
t_{1} & 0 & 0 \\ 0 & \ddots & 0 \\ 0 & 0 & t_{g} 
\end{array} \right) 
\mapsto t_{1}^{\kappa_{1}} \cdots t_{g}^{\kappa_{g}} 
$$
for $(\kappa_{1},..., \kappa_{g}) \in {\mathbb Z}^{g}$. 
Then 
$$
X^{+}(T_{g}) = \left\{ (\kappa_{1},..., \kappa_{g}) \in {\mathbb Z}^{g} \ | \ 
\kappa_{1} \geq \cdots \geq \kappa_{g} \geq 0 \right\} 
$$
becomes the set of dominant weights with respect to $B_{2g}$. 
Let $\kappa$ be an element of $X^{+}(T_{g})$ 
which is naturally regarded as regular functions on $B_{g}$ and on $B_{2g}$. 
Then    
\begin{eqnarray*}
W_{\kappa} & := & 
{\rm Ind}_{B_{g}}^{GL_{g}} (- \kappa) = 
\left\{ \left. \phi \in \Gamma \left( {\cal O}_{GL_{g}} \right) \ \right| \ 
\phi(ab) = \kappa(b) \phi(a) \ (b \in B_{g}) \right\}, 
\\
V_{\kappa} & := & 
{\rm Ind}_{B_{2g}}^{Sp_{2g}} (- \kappa) = 
\left\{ \left. \psi \in \Gamma \left( {\cal O}_{Sp_{2g}} \right) \ \right| \ 
\psi(ab) = \kappa(b) \psi(a) \ (b \in B_{2g}) \right\}   
\end{eqnarray*}
are representation spaces of $GL_{g}$, $Sp_{2g}$ by 
\begin{eqnarray*}
\phi(a) & \mapsto & (\alpha \cdot \phi)(a) = \phi (\alpha^{-1} a) \ \ 
\left( \phi \in W_{\kappa}, \ \alpha \in GL_{g} \right), 
\\
\psi(a) & \mapsto & (\alpha \cdot \psi)(a) = \psi (\alpha^{-1} a) \ \ 
\left( \psi \in V_{\kappa}, \ \alpha \in Sp_{2g} \right)
\end{eqnarray*}
respectively. 
The duals $W_{\kappa}^{*}$ (resp. $V_{\kappa}^{*}$) of $W_{\kappa}$ 
(resp. $V_{\kappa}$) are called the {\it universal representations} 
of highest weight  $\kappa$ (cf. [16, 5.1.3 and 8.1.2]), 
and hence the highest weight of $W_{\kappa}$ (resp. $V_{\kappa}$) are 
$(-\kappa_{g},..., -\kappa_{1})$ (resp. $\kappa$). 
By construction, 
$W_{\kappa}$, $V_{\kappa}$ give rational homomorphisms of $GL_{g}$, $Sp_{2g}$ respectively over any base ring, 
and for each $h \in {\mathbb Z}$, 
$W_{\kappa - h(1,..., 1)} \cong W_{\kappa} \otimes \det^{\otimes h}$. 
Over a field of characteristic $0$, 
$W_{\kappa}^{*}$ (resp. $V_{\kappa}^{*}$) are realized as direct summands of 
certain tensor products of $W$ (resp. $V$) associated with $\kappa$, 
and hence $W_{\kappa}$ can be regarded as a direct summand of 
$V_{\kappa}^{*} \cong V_{\kappa}$. 

If a linear map $\pi : V \rightarrow W$ satisfies that  
$W \hookrightarrow V \stackrel{\pi}{\rightarrow} W$ is the identity map on $W$ 
and that ${\rm Ker}(\pi)$ is anisotropic for the symplectic form on $V$, 
then $\pi$ gives a decomposition $V = W \oplus {\rm Ker}(\pi)$ 
compatible with the symplectic form. 
This decomposition induces an inclusion $GL_{g} \hookrightarrow Sp_{2g}$, 
and hence by the associated pullback, 
one has a ring homomorphism   
$\Gamma \left( {\cal O}_{Sp_{2g}} \right) \rightarrow 
\Gamma \left( {\cal O}_{GL_{g}} \right)$ 
which gives a $GL_{g}$-equivariant map $V_{\kappa} \rightarrow W_{\kappa}$. 
\vspace{2ex}

{\bf 2.2. Modular variety.} \ 
We review results of Chai and Faltings \cite{FC} 
on the moduli space of abelian varieties and its compactifications. 
For positive integers $g$ and $N$, 
let $\zeta_{N}$ be a primitive $N$th root of $1$, 
and ${\cal A}_{g,N}$ be the moduli stack classifying principally polarized 
abelian schemes of relative dimension $g$ with symplectic level $N$ structure. 
Then ${\cal A}_{g,N}$ is a smooth algebraic stack 
over ${\mathbb Z} \left[ 1/N, \zeta_{N} \right]$ 
of relative dimension $g(g+1)/2$, 
and becomes a fine moduli scheme if $N \geq 3$. 
Furthermore, the associated complex orbifold ${\cal A}_{g,N}({\mathbb C})$ 
is represented as the quotient space ${\cal H}_{g}/\Gamma(N)$ 
of the Siegel upper half space ${\cal H}_{g}$ of degree $g$ 
by the integral symplectic group 
$$
\Gamma(N) = 
\left\{ \gamma = \left( \begin{array}{cc} 
A_{\gamma} & B_{\gamma} \\ C_{\gamma} & D_{\gamma} \end{array} \right) 
\in Sp_{2g}({\mathbb Z}) 
\ \left| \ \begin{array}{ll} 
A_{\gamma} \equiv D_{\gamma} \equiv 1_{g} & {\rm mod}(N) \\ 
B_{\gamma} \equiv C_{\gamma} \equiv 0 & {\rm mod}(N) 
\end{array} \right. \right\} 
$$ 
of degree $g$ and level $N$ which acts on ${\cal H}_{g}$ as 
$$
{\cal H}_{g} \ni Z \mapsto \gamma(Z) = 
\left( A_{\gamma} Z + B_{\gamma} \right)
\left( C_{\gamma} Z + D_{\gamma} \right)^{-1} \in {\cal H}_{g} 
\ \left( \gamma \in \Gamma(N) \right). 
$$
Let $\pi : {\cal X} \rightarrow {\cal A}_{g,N}$ be the universal abelian scheme 
with $0$-section $s$, 
denote by ${\mathbb E}$ the {\it Hodge bundle} of rank $g$ defined as 
$\pi_{*} \left( \Omega^{1}_{{\cal X}/{\cal A}_{g,N}} \right) = 
s^{\ast} \left( \Omega^{1}_{{\cal X}/{\cal A}_{g,N}} \right)$, 
and by $\omega = \det({\mathbb E})$ the {\it Hodge line bundle}. 

For a smooth and $GL \left( {\mathbb Z}^{g} \right)$-admissible 
polyhedral cone decomposition of the space of positive semi-definite 
symmetric bilinear forms on ${\mathbb R}^{g}$, 
Chai and Faltings \cite[Chapter IV]{FC} construct the associated 
smooth compactification $\overline{\cal A}_{g,N}$ of ${\cal A}_{g,N}$, 
and the semi-abelian scheme ${\cal G}$ with $0$-section $s$ over 
$\overline{\cal A}_{g,N}$ extending ${\cal X} \rightarrow {\cal A}_{g,N}$. 
Then 
$\overline{\omega} = \det \left( s^{\ast} 
\left( \Omega^{1}_{{\cal G}/\overline{\cal A}_{g,N}} \right) \right)$ 
is an extension of $\omega = \det \left( {\mathbb E} \right)$ to 
$\overline{\cal A}_{g,N}$, and 
$$
{\cal A}_{g,N}^{\ast} = 
{\rm Proj}\left( \bigoplus_{h \geq 0} H^{0} \left( 
\overline{\cal A}_{g,N}, \overline{\omega}^{\otimes h} \right) \right)
$$
is a projective scheme over ${\mathbb Z} \left[ 1/N, \zeta_{N} \right]$ 
called {\it Satake's minimal compactification}. 
It is shown in \cite[Chapter IV, 6.8]{FC} that any geometric fiber of 
$\overline{\cal A}_{g,N}$ is irreducible, 
and hence ${\cal A}_{g,N}$ has the same property. 

Assume that $N \geq 3$. 
Then ${\cal A}_{g,N}^{\ast}$ contains ${\cal A}_{g,N}$, 
and its complement has a natural stratification by locally closed subschemes, 
each of which is isomorphic to ${\cal A}_{i,N}$ $(0 \leq i \leq g-1)$. 
Therefore, the relative codimension 
$$
{\rm codim}_{{\mathbb Z} \left[ 1/N, \zeta_{N} \right]} 
\left( {\cal A}_{g,N}^{\ast} - {\cal A}_{g,N}, \ {\cal A}_{g,N}^{\ast} \right) 
$$ 
over ${\mathbb Z} \left[ 1/N, \zeta_{N} \right]$ of 
${\cal A}_{g,N}^{\ast} - {\cal A}_{g,N}$ 
in ${\cal A}_{g,N}^{\ast}$ becomes 
$$
\frac{g(g+1)}{2} - \frac{(g-1)g}{2} = g
$$
which is greater than $1$ if $g > 1$. 
Furthermore, there is a natural morphism 
$\overline{\cal A}_{g,N} \rightarrow {\cal A}^{\ast}_{g,N}$ 
(which is an isomorphism if $g = 1$) 
extending the identity map on ${\cal A}_{g,N}$ 
such that $\overline{\omega}$ is the pullback by this morphism 
of the tautological line bundle $\omega^{\ast}$ on ${\cal A}^{\ast}_{g,N}$. 
\vspace{2ex}

{\bf 2.3. CM point.} \ 
Let $\varphi : S \rightarrow {\cal A}_{g,N}$ be a morphism of schemes 
over ${\mathbb Z} \left[ 1/N, \zeta_{N} \right]$ 
which becomes a $R$-rational point on ${\cal A}_{g,N}$ if $S = {\rm Spec}(R)$ 
for a ${\mathbb Z} \left[ 1/N, \zeta_{N} \right]$-algebra $R$. 
Then as the associated object, 
there is an abelian scheme $X$ over $S$ with principal polarization $\lambda$ 
and symplectic level $N$ structure $\sigma$. 
A {\it test object} (resp. an {\it extended test object}) over $S$ 
associated with a morphism $\varphi : S \rightarrow {\cal A}_{g,N}$ is 
the above $(X, \lambda, \sigma)$ together with basis of regular $1$-forms 
on $X/S$ (resp. basis of $H_{\rm DR}^{1}(X/S)$). 
By definition, 
any element of ${\cal M}_{\rho}(R)$ is evaluated as an element of 
${\mathcal O}_{S}^{d}$ at each test object over an $R$-scheme $S$, 
where this evaluation is functorial on $S$ and equivariant for $\rho$ 
under base changes of regular $1$-forms. 

For a field extension $k$ of ${\mathbb Q}(\zeta_{N})$, 
a $k$-rational point $\alpha$ on ${\cal A}_{g,N}$ corresponding to a CM abelian variety  $X$ is called a {\it CM point over $k$} if the following conditions hold: 
\begin{itemize}

\item 
The endomorphism ${\mathbb Q}$-algebra ${\rm End}_{k}(X) \otimes {\mathbb Q}$ 
is isomorphic to the direct sum $\bigoplus_{i} L_{i}$, 
where $L_{i}$ are CM fields, i.e., 
totally imaginary quadratic extensions of totally real number fields $K_{i}$. 

\item 
There are algebra homomorphisms 
$\varphi_{i} : L_{i} \otimes k \rightarrow K_{i} \otimes k$ such that 
$$
x \otimes y \mapsto \left( \varphi_{i}(x \otimes y), 
\varphi_{i} \left( \iota_{i}(x) \otimes y \right) \right) 
\ \left( x \in L_{i}, y \in k \right) 
$$ 
give rise to isomorphisms 
$L_{i} \otimes k \stackrel{\sim}{\rightarrow} 
K_{i} \otimes k \oplus K_{i} \otimes k$, 
where $\iota_{i}$ denotes the involution of $L_{i}$ over $K_{i}$. 

\end{itemize}

Note that any CM abelian variety can be defined over a number field, 
and has potentially good reduction at all finite places. 
Therefore, 
for any CM point $\alpha$ on ${\cal A}_{g,N}$ and any rational prime $p$, 
there is an (extended) test object $\widetilde{\alpha}$ associated with $\alpha$ 
over an algebra which is a finite ${\mathbb Z}_{(p)}$-module, 
where ${\mathbb Z}_{(p)}$ denotes the valuation ring of ${\mathbb Q}$ at $p$. 
\vspace{2ex}

{\bf 2.4. Modular forms.} \ 
In what follows, we assume that 
$$
g > 1, \ N \geq 3. 
$$

First, following \cite[2.2.1]{G1} 
we give the process of twisting a locally free sheaf by 
a linear representation. 
Let $X$ be a scheme, and ${\cal F}$ be a locally free sheaf on $X$ 
of rank $n$. 
Take $\left\{ U_{i} \right\}_{i \in I}$ be an open cover of $X$ 
trivializing ${\cal F}$. 
Then the natural isomorphism 
${\cal F}|_{U_{i} \cap U_{j}} \cong {\cal F}|_{U_{j} \cap U_{i}}$ 
gives rise to the transition function 
$g_{ij} \in GL_{n} \left( {\cal O}_{X}|_{U_{i} \cap U_{j}} \right)$ 
satisfying the cocycle condition. 
Let $\rho : GL_{n} \rightarrow GL_{m}$ be a rational homomorphism over 
a ${\mathbb Z}$-algebra $R$. 
Then we construct a locally free ${\cal O}_{X} \otimes R$-module 
${\cal F}_{\rho}$ on $X \otimes R$ as 
${\cal F}_{\rho}|_{U_{i}} = 
\left( \left( {\cal O}_{X} \otimes R \right)|_{U_{i}} \right)^{m}$, 
where the isomorphism 
${\cal F_{\rho}}|_{U_{i} \cap U_{j}} \cong 
{\cal F_{\rho}}|_{U_{j} \cap U_{i}}$ 
is given by 
$\rho(g_{ij}) \in GL_{m} 
\left( \left( {\cal O}_{X} \otimes R\right) |_{U_{i} \cap U_{j}} \right)$. 

For a ${\mathbb Z} \left[ 1/N, \zeta_{N} \right]$-algebra $R$, 
a positive integer $d$ and a rational homomorphism 
$\rho : GL_{g} \rightarrow GL_{d}$ over $R$, 
let ${\mathbb E}_{\rho}$ be the locally free sheaf on 
$$
{\cal A}_{g,N} \otimes R = 
{\cal A}_{g,N} \otimes_{{\mathbb Z} \left[ 1/N, \zeta_{N} \right]} R
$$ 
obtained from twisting the Hodge bundle ${\mathbb E}$ by $\rho$. 
If $\rho$ is obtained from $\kappa \in {\mathbb Z}^{g}$, 
then we put ${\mathbb E}_{\kappa} = {\mathbb E}_{\rho}$, 
and denote this rank $d$ by $d({\mathbb E}_{\kappa})$. 
\vspace{2ex}

{\sc Definition 2.1.} 
Let $R$ be a ${\mathbb Z} \left[ 1/N, \zeta_{N} \right]$-algebra. 
For a rational homomorphism $\rho : GL_{g} \rightarrow GL_{d}$ over $R$, 
we put 
$$
{\cal M}_{\rho}(R) = 
H^{0} \left( {\cal A}_{g,N} \otimes R, {\mathbb E}_{\rho} \right), 
$$
and call these elements {\it Siegel modular forms over $R$ of weight $\rho$} 
(and degree $g$, level $N$). 
If $\rho = \omega^{\otimes h} : GL_{g} \rightarrow {\mathbb G}_{m}$, 
then we put ${\cal M}_{h}(R) = {\cal M}_{\omega^{\otimes h}}(R)$, 
and call these elements of weight $h$. 
More generally, for an $R$-module $M$, 
the space of {\it Siegel modular forms with coefficients in $M$ of weight $\rho$} is defined as 
$$
{\cal M}_{\rho}(M) = H^{0} \left( {\cal A}_{g,N} \otimes R, 
{\mathbb E}_{\rho} \otimes_{R} M \right). 
$$

We consider the case where $R = {\mathbb C}$. 
For each $Z \in {\cal H}_{g}$, 
let 
$$
{\cal X}_{Z} = {\mathbb C}^{g} / 
({\mathbb Z}^{g} + {\mathbb Z}^{g} \cdot Z)
$$ 
be the corresponding abelian variety over ${\mathbb C}$, 
and $(u_{1}, ..., u_{g})$ be the natural coordinates on the universal cover 
${\mathbb C}^{g}$ of ${\cal X}_{Z}$. 
Then ${\mathbb E}$ is trivialized over ${\cal H}_{g}$ by 
$du_{1}, ..., du_{g}$. 
For 
${\displaystyle \gamma = \left( \begin{array}{cc} 
A_{\gamma} & B_{\gamma} \\ C_{\gamma} & D_{\gamma} \end{array} 
\right) \in \Gamma(N)}$, 
$$
{\cal X}_{Z} \stackrel{\sim}{\rightarrow} {\cal X}_{\gamma(Z)} ; \ 
\mbox{}^{t}(u_{1}, ..., u_{g}) \mapsto 
\left( C_{\gamma} Z + D_{\gamma} \right)^{-1} \cdot 
\mbox{}^{t}(u_{1}, ..., u_{g}), 
$$
and hence $\gamma$ acts equivariantly on the trivialization of ${\mathbb E}$ 
over ${\cal H}_{g}$ as the left multiplication by 
$\left( C_{\gamma} Z + D_{\gamma} \right)^{-1}$. 
Therefore, $\gamma$ acts equivariantly on the induced trivialization of 
${\mathbb E}_{\rho}$ over ${\cal H}_{g}$ as the left multiplication by 
$\rho \left( C_{\gamma} Z + D_{\gamma} \right)^{-1}$. 
Then $f \in {\cal M}_{\rho}({\mathbb C})$ is a complex analytic section of 
${\mathbb E}_{\rho}$ on 
${\cal A}_{g,N}({\mathbb C}) = {\cal H}_{g}/\Gamma(N)$, 
and hence is a ${\mathbb C}^{d}$-valued holomorphic function on ${\cal H}_{g}$ 
satisfying the $\rho$-automorphic condition: 
$$
f(Z) = \rho \left( C_{\gamma} Z + D_{\gamma} \right)^{-1} \cdot 
f \left( \gamma(Z) \right) \ 
\left( Z \in {\cal H}_{g}, \ 
\gamma = \left( \begin{array}{cc} 
A_{\gamma} & B_{\gamma} \\ C_{\gamma} & D_{\gamma} \end{array} \right) \in 
\Gamma(N) \right) 
$$
which is equivalent to that 
$f \left( \gamma(Z) \right) = 
\rho \left( C_{\gamma} Z + D_{\gamma} \right) \cdot f(Z)$. 
Furthermore, the value of $f$ at a test object 
$\left( X, \lambda, \alpha; w_{1},..., w_{g} \right)$ over a subfield $k$ of ${\mathbb C}$ becomes $\rho(G) \cdot f(Z) \in k^{d}$, 
where 
$\mbox{}^{t}( du_{1},..., du_{g}) = G \cdot \mbox{}^{t}( w_{1},..., w_{g})$. 

Let $\iota : {\cal A}_{g,N} \hookrightarrow {\cal A}_{g,N}^{\ast}$ 
be the natural inclusion, 
and let ${\mathbb E}_{\rho}^{\ast}$ be the direct image (or pushforward) 
$\iota_{\ast} \left( {\mathbb E}_{\rho} \right)$ which is defined as 
a sheaf on ${\cal A}_{g,N}^{\ast} \otimes R$ satisfying that 
${\mathbb E}_{\rho}^{\ast}(U) = 
{\mathbb E}_{\rho} \left( \iota^{-1}(U) \right)$ 
for open subsets $U$ of ${\cal A}_{g,N}^{\ast} \otimes R$. 
This implies immediately that 
$$
{\cal M}_{\rho}(R) = 
H^{0} \left( {\cal A}_{g,N}^{\ast} \otimes R, 
{\mathbb E}_{\rho}^{\ast} \right). 
$$
Furthermore, based on that 
${\rm codim}_{{\mathbb Z} \left[ 1/N, \zeta_{N} \right]} 
\left( {\cal A}_{g,N}^{\ast} - {\cal A}_{g,N}, \ {\cal A}_{g,N}^{\ast} \right) 
> 1$, 
Ghitza \cite[Theorem 3]{G2} proved that ${\mathbb E}_{\rho}^{\ast}$ is 
a coherent sheaf on ${\cal A}_{g,N}^{\ast} \otimes R$. 
From this fact, it is shown in \cite[Theorem 1]{I} that 
${\cal M}_{\rho}(R)$ is a finitely generated $R$-module, 
and that ${\cal M}_{\rho}({\mathbb C})$ consists of 
${\mathbb C}^{d}$-valued holomorphic functions on ${\cal H}_{g}$ 
satisfying the $\rho$-automorphic condition. 
\vspace{2ex}

{\bf 2.5. Fourier expansion.} \ 
Let $q_{ij}$ $ \left( 1 \leq i,j \leq g \right)$ be variables with symmetry 
$q_{ij} = q_{ji}$. 
Then in \cite{Mu}, 
Mumford constructs a semi-abelian scheme formally represented as 
$$
\left. \left( {\mathbb G}_{m} \right)^{g} \right/ 
\langle \left. (q_{ij})_{1 \leq i \leq g} \ \right| \ 1 \leq j \leq g \rangle 
; \ \left( {\mathbb G}_{m} \right)^{g} = 
{\rm Spec} \left( {\mathbb Z} \left[ x_{1}^{\pm 1},...,x_{g}^{\pm 1} \right] 
\right) 
$$
over 
$$
{\mathbb Z} \left[ q_{ij}^{\pm 1} \ (i \neq j) \right] 
\left[ \left[ q_{11},...,q_{gg} \right] \right]. 
$$
This becomes an abelian scheme which is called {\it Mumford's abelian scheme} 
over  
$$
{\mathbb Z} \left[ q_{ij}^{\pm 1} \ (i \neq j) \right] 
\left[ \left[ q_{11},...,q_{gg} \right] \right] 
\left[ 1/q_{11},...,1/q_{gg} \right] 
$$
with principal polarization corresponding to the multiplicative form 
$$
\left( (a_{1},...,a_{g}), (b_{1},...,b_{g}) \right) \mapsto 
\prod_{1 \leq i,j \leq g} q_{ij}^{a_{i}b_{j}} 
$$
on ${\mathbb Z}^{g} \times {\mathbb Z}^{g}$. 
Hence for each $0$-dimensional cusp $c$ on ${\cal A}_{g,N}^{\ast}$, 
this polarized abelian scheme over 
$$
{\cal R}_{g,N} =  
{\mathbb Z} \left[ 1/N, \ \zeta_{N}, \ q_{ij}^{\pm 1/N} \ (i \neq j) \right] 
\left[ \left[ q_{11}^{1/N},...,q_{gg}^{1/N} \right] \right] 
\left[ 1/q_{11},...,1/q_{gg} \right] 
$$
has the associated symplectic level $N$ structure, 
and $\omega_{i} = d x_{i}/x_{i}$ $(1 \leq i \leq g)$ 
form a basis of regular $1$-forms. 
Taking the pullback by the associated morphism 
${\rm Spec} \left( {\cal R}_{g,N} \right) \rightarrow {\cal A}_{g,N}$, 
${\mathbb E}$ is trivialized by the basis $\omega_{1},...,\omega_{g}$, 
and hence ${\mathbb E}_{\rho}$ is also trivialized over 
${\rm Spec} \left( {\cal R}_{g,N} \otimes R \right) = 
{\rm Spec} \left( {\cal R}_{g,N} 
\otimes_{{\mathbb Z}[1/N, \zeta_{N}]} R \right)$. 
In what follows, we fix such a trivialization: 
$$
{\mathbb E}_{\rho} \times_{{\cal A}_{g,N} \otimes R} 
{\rm Spec} \left( {\cal R}_{g,N} \otimes R \right) = 
\left( {\cal R}_{g,N} \otimes R \right)^{d}. 
$$
Then for an $R$-module $M$, 
the evaluation on Mumford's abelian scheme gives a homomorphism 
$$
F_{c} : {\cal M}_{\rho}(M) \rightarrow 
\left( {\cal R}_{g,N} \otimes_{{\mathbb Z}[1/N, \zeta_{N}]} M \right)^{d}
$$
which we call the {\it Fourier expansion map} associated with $c$. 
Furthermore, it is shown in \cite[Theorem 2]{I} that $F_{c}$ satisfies 
the following $q$-expansion principle: 
\vspace{2ex}

\begin{it} 
If $M'$ is an $R$-submodule of $M$ and $f \in {\cal M}_{\rho}(M)$ 
satisfies that 
$$
F_{c}(f) \in 
\left( {\cal R}_{g,N} \otimes_{{\mathbb Z}[1/N, \zeta_{N}]} M' \right)^{d}, 
$$
then $f \in {\cal M}_{\rho}(M')$. 
\end{it}
\vspace{2ex}

\noindent
This result was already shown by Harris \cite[4.8, Theorem]{H} 
in the case where $M$ is a field extension of a field $M'$ 
containing ${\mathbb Q}(\zeta_{N})$. 

Assume that $M = {\mathbb C}$ and $c$ is associated with $\sqrt{-1} \infty$. 
Then by the substitution $q_{ij} = \exp \left( 2 \pi \sqrt{-1} z_{ij} \right)$ 
for $Z = (z_{ij})_{i,j} \in {\cal H}_{g}$, 
Mumford's abelian scheme becomes ${\cal X}_{Z}$, 
and hence $F_{c}$ becomes the analytic Fourier expansion map times 
$\rho \left( 2 \pi \sqrt{-1} \cdot 1_{g} \right)$. 
Since each $f(Z) \in {\cal M}_{\rho}({\mathbb C})$ is 
a ${\mathbb C}^{d}$-valued holomorphic function of $Z \in {\cal H}_{g}$ 
and is invariant under $Z \mapsto Z + N \cdot I$ 
for any integral and symmetric $g \times g$ matrix $I$, 
and hence 
$$
F_{c}(f) = \sum_{T} a(T) \cdot 
\exp \left( 2 \pi \sqrt{-1} {\rm tr}(TZ)/N \right) = 
\sum_{T} a(T) \cdot \mbox{\boldmath $q$}^{T/N} \ 
\left( a(T) \in {\mathbb C}^{d} \right), 
$$
where $T = \left( t_{ij} \right)_{i,j}$ runs over half-integral 
symmetric $g \times g$ matrices, 
and 
$$
\mbox{\boldmath $q$}^{T/N} = 
\prod_{1 \leq i < j \leq g} \left( q_{ij}^{1/N} \right)^{2 t_{ij}} 
\prod_{1 \leq i \leq g} \left( q_{ii}^{1/N} \right)^{t_{ii}}. 
$$
Furthermore, as is shown in the Cartan Seminar 4-04, 
$a(T) = 0$ if $T$ is not positive semi-definite. 
\vspace{2ex}

{\bf 2.6. Nearly modular forms.} \ 
Let ${\cal H}_{\rm DR}^{1} \left( {\cal X} / {\cal A}_{g,N} \right)$ be the sheaf 
of de Rham cohomology groups of ${\cal X}/{\cal A}_{g,N}$, 
and define the {\it de Rham bundle} as 
$$
{\mathbb D} = {\cal R}_{\rm DR}^{1} \pi \left( {\cal X} / {\cal A}_{g,N} \right) = 
\pi_{*} \left( {\cal H}_{\rm DR}^{1} \left( {\cal X} / {\cal A}_{g,N} \right) \right)
$$
which is a locally free sheaf on ${\cal A}_{g,N}$ of rank $2g$ 
with canonical symplectic form. 
Then one has a canonical exact sequence 
$$
0 \rightarrow {\mathbb E} \rightarrow {\mathbb D} \rightarrow {\mathbb D}/{\mathbb E} 
\rightarrow 0, 
$$
and the quotient ${\mathbb D}/{\mathbb E}$ is locally free of rank $g$. 
The Gauss-Manin connection 
$$
\nabla : {\mathbb D} \rightarrow {\mathbb D} \otimes \Omega_{{\cal A}_{g,N}}
$$
defines 
${\cal T}_{{\cal A}_{g,N}} \rightarrow {\rm End}_{{\cal O}_{{\cal A}_{g,N}}} ({\mathbb D})$ 
which, 
together with the above exact sequence, 
gives the Kodaira-Spencer isomorphism 
$$
{\cal T}_{{\cal A}_{g,N}} \stackrel{\sim}{\rightarrow} 
{\rm Hom}_{{\cal O}_{{\cal A}_{g,N}}} \left( {\mathbb E}, {\mathbb D}/{\mathbb E} \right).
$$

Let $\kappa$ be an element of $X^{+}(T_{g})$, 
and denote by $V_{\kappa}$ the universal representation of highest weight $\kappa$. 
Then one can obtain the associated locally free sheaf ${\mathbb D}_{\kappa}$ 
on ${\cal A}_{g,N}$ whose rank is denoted by $d({\mathbb D}_{\kappa})$. 
Furthermore, for $h \in {\mathbb Z}$, put 
$$
{\mathbb D}_{(\kappa, h)} = {\mathbb D}_{\kappa} \otimes \det({\mathbb E})^{\otimes h} 
$$
which is also a locally free sheaf on ${\cal A}_{g,N}$ 
with rank $d({\mathbb D}_{\kappa})$. 
\vspace{2ex}

{\sc Definition 2.2.} 
Let $R$ be a ${\mathbb Z} \left[ 1/N, \zeta_{N} \right]$-algebra. 
Then for $\kappa \in X^{+}(T_{g})$ and $h \in {\mathbb Z}$, 
we put 
$$
{\cal N}_{(\kappa, h)}(R) = 
H^{0} \left( {\cal A}_{g,N} \otimes R, {\mathbb D}_{(\kappa, h)} \right), 
$$
and call these elements {\it nearly Siegel modular forms over $R$ 
of weight $(\kappa, h)$} (and degree $g$, level $N$). 
More generally, for an $R$-module $M$,  
we call 
$$ 
{\cal N}_{(\kappa, h)}(M) = 
H^{0} \left( {\cal A}_{g,N} \otimes R, {\mathbb D}_{(\kappa, h)} \otimes_{R} M \right). 
$$ 
the space of {\it nearly Siegel modular forms with coefficients in $M$ 
of weight $(\kappa, h)$}. 
\vspace{2ex}  

{\sc Theorem 2.3.} 
\begin{it} 
The $R$-module ${\cal N}_{(\kappa, h)}(R)$ is finitely generated. 
\end{it} 
\vspace{2ex} 

{\it Proof.} 
Let $\iota : {\cal A}_{g,N} \hookrightarrow {\cal A}_{g,N}^{*}$ be the natural inclusion. 
Then by \cite[Theorem 3]{G2},  
$\iota_{*} \left( {\mathbb D}_{(\kappa, h)} \right)$ 
is a coherent sheaf on ${\cal A}_{g,N}$, 
and hence 
$$
{\cal N}_{(\kappa, h)}(R) = 
H^{0} \left( {\cal A}_{g,N}^{*} \otimes R, 
\iota_{*} \left( {\mathbb D}_{(\kappa, h)} \right) \right)
$$
is a finitely generated $R$-module. 
\ $\square$ 
\vspace{2ex} 

As in 2.5, let  $\{ \omega_{i} \ | \ 1 \leq i \leq g \}$ be the canonical basis 
of the Mumford's abelian scheme. 
Then there exist $\eta_{i}$ $(1 \leq i \leq g)$ such that 
$$
\nabla(\omega_{i}) = \sum_{j=1}^{g} \frac{d q_{ij}}{q_{ij}} \eta_{j}, 
$$ 
and $\{ \omega_{i}, \eta_{i} \ | \ 1 \leq i \leq g \}$ gives a basis of ${\mathbb D}$ 
over ${\cal A}_{g,N}$. 
By using this basis, 
one has a trivialization of ${\mathbb D}_{\kappa}$ over ${\cal A}_{g,N}$ 
and that of $\det({\mathbb E})$ by $\omega_{1} \wedge \cdots \wedge \omega_{g}$. 
Therefore, there exists the Fourier expansion map  
$$
{\cal F}_{c} : {\cal N}_{(\kappa, h)}(M) \rightarrow 
\left( {\cal R}_{g,N} 
\otimes_{{\mathbb Z}[1/N, \zeta_{N}]} M \right)^{d({\mathbb D}_{\kappa})} 
$$ 
which is obtained as the evaluation map on the Mumford's abelian scheme. 
\vspace{2ex}

{\sc Theorem 2.4.} 
\begin{it} 
The Fourier expansion map ${\cal F}_{c}$ satisfies 
the following $q$-expansion principle: 
If $M'$ is an $R$-submodule of $M$ and $f \in {\cal N}_{(\kappa, h)}(M)$ satisfies 
$$
{\cal F}_{c}(f) \in \left( {\cal R}_{g,N} \otimes_{{\mathbb Z}[1/N, 
\zeta_{N}]} M' \right)^{d({\mathbb D}_{\kappa})}, 
$$
then $f \in {\cal N}_{(\kappa, h)}(M')$. 
\end{it}
\vspace{2ex}

{\it Proof.} 
The proof goes in the same way to that of \cite[Theorem 2]{I}. 
Since ${\cal A}_{g,N}$ is smooth over ${\mathbb Z} \left[ 1/N, \zeta_{N} \right]$, 
${\mathbb D}_{(\kappa, h)} \otimes R$ is flat over $R$, 
and hence ${\cal N}_{\kappa}$ is left-exact in $M$. 
Therefore, to prove the assertion, 
it is enough to show the injectivity of ${\cal F}_{c}$. 
Let $P$ be a point on ${\cal A}_{g,N}$, 
and ${\rm Spec}(A)$ be an open neighborhood of $P$ in ${\cal A}_{g,N}$. 
Then by the properties of ${\cal A}_{g,N}$ mentioned above, 
${\rm Spec}(A)$ is smooth over 
${\mathbb Z} \left[ 1/N, \zeta_{N} \right]$, 
and its geometric fibers are all irreducible. 
Since the morphism 
${\rm Spec} \left( {\cal R}_{g,N} \right) \rightarrow {\cal A}_{g,N}$ 
associated with Mumford's abelian scheme is dominant 
over any geometric fiber of 
${\rm Spec} \left( {\mathbb Z} \left[ 1/N, \zeta_{N} \right] \right)$, 
the induced morphism 
$$
{\rm Spec} \left( {\cal R}_{g,N} \right) \times_{{\cal A}_{g,N}} {\rm Spec}(A) 
\rightarrow {\rm Spec}(A) 
$$
also has the same property. 
Put 
$$
{\rm Spec} \left( {\cal R}_{g,N/A} \right) = 
{\rm Spec} \left( {\cal R}_{g,N} \right) \times_{{\cal A}_{g,N}} {\rm Spec}(A). 
$$
Then by the Chinese remainder theorem, 
for any ideal $I$ of ${\mathbb Z} \left[ 1/N, \zeta_{N} \right]$, 
$$
A \otimes \left( {\mathbb Z} \left[ 1/N, \zeta_{N} \right] / I \right) 
\rightarrow 
{\cal R}_{g,N/A} \otimes 
\left( {\mathbb Z} \left[ 1/N, \zeta_{N} \right] / I \right) 
$$
is injective. 
Since $A$ and ${\cal R}_{g,N/A}$ are flat 
${\mathbb Z} \left[ 1/N, \zeta_{N} \right]$-modules, 
if $M$ is a finitely generated 
${\mathbb Z} \left[ 1/N, \zeta_{N} \right]$-module, 
then $A \otimes M \rightarrow {\cal R}_{g,N/A} \otimes M$ is injective. 
Hence this property holds for any 
${\mathbb Z} \left[ 1/N, \zeta_{N} \right]$-module $M$ 
because any ${\mathbb Z} \left[ 1/N, \zeta_{N} \right]$-module is 
the inductive limit of  finitely generated 
${\mathbb Z} \left[ 1/N, \zeta_{N} \right]$-modules, 
and tensor product commutes with inductive limits. 
We may assume that ${\mathbb D}_{(\kappa, h)}$ is trivialized 
on ${\rm Spec}(A \otimes R)$, 
and hence 
$$
\left( {\mathbb D}_{(\kappa, h)} |_{{\rm Spec}(A \otimes R)} \right)  \otimes_{R} M 
\rightarrow \left( {\cal R}_{g,N/A} \otimes M \right)^{d({\mathbb D}_{\kappa})} 
$$
is injective for any $R$-module $L$. 
Therefore, 
any $f \in {\cal N}_{(\kappa, h)}(M)$ satisfying that ${\cal F}_{c}(f) = 0$ 
vanishes around any point $P$ on ${\cal A}_{g,N}$, 
and hence $f = 0$. \ $\square$    
\vspace{2ex}

{\sc Theorem 2.5.} 
\begin{it} 
Let $R$ be a ${\mathbb Z} \left[ 1/N, \zeta_{N} \right]$-algebra, 
and $f$ be an element of ${\cal N}_{(\kappa, h)}(R)$. 
Then for an extended test object $\widetilde{\alpha}$ over $R$ associated with 
a point $\alpha$ on ${\cal A}_{g,N}$, 
the evaluation $f \left( \widetilde{\alpha} \right)$ of $f$ 
at $\widetilde{\alpha}$ belongs to $R^{d({\mathbb D}_{\kappa})}$. 
\end{it}
\vspace{2ex}

{\it Proof.} 
This assertion directly follows from the definition of nearly Siegel modular forms 
and their rationality. \ $\square$  
\vspace{3ex}

{\bf 3. Arithmeticity in the analytic case} 
\vspace{2ex}

{\bf 3.1. Differential operator.} \ 
First, we recall Shimura's differential operator. 
Let $R$ be a ${\mathbb Q}$-algebra, 
and identify the $2$-fold symmetric tensor product ${\rm Sym}^{2}(R^{g})$ of 
$R^{g}$ with the $R$-module of all symmetric $g \times g$ matrices with entries in $R$. 
For a positive integer $e$, 
let $S_{e} \left( {\rm Sym}^{2}(R^{g}), R^{d} \right)$ be the $R$-module of 
all polynomial maps of ${\rm Sym}^{2}(R^{g})$ into $R^{d}$ 
homogeneous of degree $e$. 
For a rational homomorphism $\rho : GL_{g} \rightarrow GL_{d}$, 
let $\rho \otimes \tau^{e}$ and $\rho \otimes \sigma^{e}$ 
be the rational homomorphisms over $R$ given by 
$$
GL_{g}(R) = {\rm Aut}_{R} (R^{g}) \rightarrow {\rm Aut}_{R} 
\left( S_{e} \left( {\rm Sym}^{2}(R^{g}), R^{d} \right) \right) 
$$ 
which are defines as 
$$
\left[ \left( \rho \otimes \tau^{e} \right)(\alpha)(h) \right](u) = 
\rho(\alpha) h \left( \mbox{}^{t} \alpha \cdot u \cdot \alpha \right)
$$
and
$$
\left[ \left( \rho \otimes \sigma^{e} \right)(\alpha)(h) \right](u) = 
\rho(\alpha) h \left( \alpha^{-1} \cdot u \cdot \mbox{}^{t} \alpha^{-1} \right)
$$
respectively for 
$\alpha \in GL_{g}(R)$, 
$h \in S_{e} \left( {\rm Sym}^{2}(R^{g}), R^{d} \right)$, 
$u \in {\rm Sym}^{2}(R^{g})$. 
In particular, for $\alpha \in GL_{g}$, 
$\tau^{e}(\alpha)$ (resp. $\sigma^{e}(\alpha)$) consists of polynomials of 
entries of $\alpha$ (resp. $\alpha^{-1}$). 
Furthermore, let 
$$
\theta^{e} : 
S_{e} \left( {\rm Sym}^{2}(R^{g}), 
S_{e} \left( {\rm Sym}^{2}(R^{g}), R^{d} \right) \right) 
\rightarrow R^{d}
$$
be the contraction map defined in \cite[14.1]{Sh} as 
$\theta^{e}(h) = \sum_{i} h \left( u_{i}, v_{i} \right)$, 
where $\{ u_{i} \}$ and $\{ v_{i} \}$ are dual basis of 
${\rm Sym}^{2} \left( R^{g} \right)$ 
for the pairing $(u, v) \mapsto {\rm tr}(uv)$, 
namely ${\rm tr}(u_{i} v_{j})$ is Kronecker's delta $\delta_{ij}$. 
Then $\theta^{e}$ is $GL_{g}$-equivariant for the representations 
$\rho \otimes \sigma^{e} \otimes \tau^{e}$ and $\rho$. 

Let $f$ be a ${\mathbb C}^{d}$-valued smooth function of 
$Z = \left( z_{ij} \right)_{i,j} = X + \sqrt{-1} Y \in {\cal H}_{g}$. 
Then following \cite[Chapter III, 12]{Sh}, 
define 
$S_{1} \left( {\rm Sym}^{2}({\mathbb C}^{g}), {\mathbb C}^{g} \right)$-valued 
smooth functions $(Df)(u)$, $(Cf)(u)$ 
$\left( u = \left( u_{ij} \right)_{i,j} \in 
{\rm Sym}^{2}({\mathbb C}^{g}) \right)$ 
of $Z \in {\cal H}_{g}$ as 
\begin{eqnarray*}
(Df)(u) 
& = & 
\sum_{1 \leq i \leq j \leq g} u_{ij} 
\frac{\partial f}{\partial (2 \pi \sqrt{-1} z_{ij})}, 
\\
(Cf)(u) 
& = & 
(Df) 
\left( \left(Z - \overline{Z} \right) u 
\left(Z - \overline{Z} \right) \right), 
\end{eqnarray*}
and define 
$S_{e} \left( {\rm Sym}^{2}({\mathbb C}^{g}), {\mathbb C}^{g} \right)$-valued 
analytic functions $C^{e}(f)$, $D_{\rho}^{e}(f)$ of $Z \in {\cal H}_{g}$ as 
\begin{eqnarray*}
C^{e}(f) 
& = & 
C \left( C^{e-1}(f) \right), 
\\
D_{\rho}^{e}(f) 
& = & 
(\rho \otimes \tau^{e}) \left(Z - \overline{Z} \right)^{-1} 
C^{e} \left( \rho \left(Z - \overline{Z} \right) f \right). 
\end{eqnarray*}
It is shown in \cite[Chapter III, 12.10]{Sh} that 
if $f$ satisfies the $\rho$-automorphic condition for $\Gamma(N)$, 
then $D_{\rho}^{e}(f)(u)$ satisfies 
the $\rho \otimes \tau^{e}$-automorphic condition. 
\vspace{2ex}

{\it Remark 1.} 
The above $D_{\rho}^{e}$ becomes $\left( 2 \pi \sqrt{-1} \right)^{-e}$ 
times Shimura's original operator given in \cite{Sh}. 
\vspace{2ex}

Let $u_{1}, ..., u_{g}$ be the standard coordinates on ${\mathbb C}^{g}$, 
and $\alpha_{i}$, $\beta_{i}$ $(1 \leq i \leq g)$ be relative $1$-forms on 
${\cal X}_{Z}$ $\left( Z = (z_{ij}) \in {\cal H}_{g} \right)$ given by 
$$
\alpha_{i} \left( \sum_{j=1}^{g} a_{j} \mbox{\boldmath $e$}_{j} + 
\sum_{j=1}^{g} b_{j} \mbox{\boldmath $z$}_{j} \right) = a_{i}, 
\ 
\beta_{i} \left( \sum_{j=1}^{g} a_{j} \mbox{\boldmath $e$}_{j} + 
\sum_{j=1}^{g} b_{j} \mbox{\boldmath $z$}_{j} \right) = b_{i} 
$$ 
for each $a_{j}, b_{j} \in {\mathbb R}$, 
where 
$\mbox{\boldmath $e$}_{j} = \left( \delta_{ij} \right)_{1 \leq i \leq g}$ and 
$\mbox{\boldmath $z$}_{j} = \left( z_{j1}, ..., z_{jg} \right)$. 
Since $\alpha_{i}$, $\beta_{i}$ have constant periods for all ${\cal X}_{Z}$, 
$\nabla(\alpha_{i}) = \nabla(\beta_{i}) = 0$. 
Furthermore, one has 
$$
d u_{i} = \alpha_{i} + \sum_{j=1}^{g} z_{ij} \beta_{j}, 
\ 
d \overline{u_{i}} = 
\alpha_{i} + \sum_{j=1}^{g} \overline{z_{ij}} \beta_{j} 
$$
which implies that 
$$
\mbox{}^{t}(du_{1}, ..., du_{g}) \equiv \left( Z - \overline{Z} \right) \cdot 
\mbox{}^{t}(\beta_{1}, ..., \beta_{g}) \ \ 
{\rm mod} \left( H^{0,1} \left( {\cal X}/{\cal H}_{g} \right) \right). 
$$
Then 
$$
\omega_{i} = d \log(x_{i}) = 2 \pi \sqrt{-1} d u_{i} \ (1 \leq i \leq g), 
$$
and hence 
$$ 
\nabla(\omega_{i}) = 2 \pi \sqrt{-1} \nabla(d u_{i}) = 
2 \pi \sqrt{-1} \sum_{j=1}^{g} d z_{ij} \cdot \beta_{j} = 
2 \pi \sqrt{-1} \sum_{j=1}^{g} \frac{d q_{ij}}{q_{ij}} \beta_{j} 
$$
which implies 
$$
\eta_{i} = \beta_{i} \ (1 \leq i \leq g). 
$$
\vspace{-1ex}
 
The following proposition was obtained by Harris \cite[Section 4]{H} substantially, 
and shown by Eischen \cite[Proposition 8.5]{E1} in the unitary modular case. 
\vspace{2ex}

{\sc Proposition 3.1.} 
\begin{it} 
Let $\pi : {\cal X} \rightarrow {\cal H}_{g}$ be the family of 
complex abelian varieties given by 
$$
\pi^{-1}(Z) = {\cal X}_{Z} = 
{\mathbb C}^{g} / ({\mathbb Z}^{g} + {\mathbb Z}^{g} \cdot Z) \ 
(Z \in {\cal H}_{g}). 
$$
Then $D_{\rho}^{e}$ is obtained from the composition 
$$
{\mathbb E}_{\rho} 
\rightarrow 
{\mathbb E}_{\rho} \otimes 
\left( \Omega_{{\cal H}_{g}}^{1} \right)^{\otimes e} 
\rightarrow 
{\mathbb E}_{\rho} \otimes 
\left( {\rm Sym}^{2} \left( \pi_{*} \left( 
\Omega_{{\cal X}/{\cal H}_{g}}^{1} \right) \right) \right)^{\otimes e}. 
$$ 
Here the first map is given by the Gauss-Manin connection 
$$
\nabla : H_{\rm DR}^{1} \left( {\cal X}/{\cal H}_{g} \right) \rightarrow 
H_{\rm DR}^{1} \left( {\cal X}/{\cal H}_{g} \right) \otimes 
\Omega_{{\cal H}_{g}}^{1} 
$$
together with the projection 
$H_{\rm DR}^{1} \left( {\cal X}/{\cal H}_{g} \right) \rightarrow 
\pi_{*} \left( \Omega_{{\cal X}/{\cal H}_{g}}^{1} \right)$ 
derived from the Hodge decomposition 
$$
H_{\rm DR}^{1} \left( {\cal X}/{\cal H}_{g} \right) 
\cong 
H^{1,0} \left( {\cal X}/{\cal H}_{g} \right) \oplus 
H^{0,1} \left( {\cal X}/{\cal H}_{g} \right) 
= 
\pi_{*} \left( \Omega_{{\cal X}/{\cal H}_{g}}^{1} \right) \oplus 
\overline{\pi_{*} \left( \Omega_{{\cal X}/{\cal H}_{g}}^{1} \right)}, 
$$ 
and the second map is given by the Kodaira-Spencer isomorphism 
\end{it}
$$
\Omega_{{\cal H}_{g}}^{1} \cong 
{\rm Sym}^{2} \left( \pi_{*} \left( 
\Omega_{{\cal X}/{\cal H}_{g}}^{1} \right) \right). 
$$

{\it Proof.} \ 
We regard a ${\mathbb C}^{d}$-valued smooth function $f$ of 
$Z \in {\cal H}_{g}$ with $\rho$-automorphic condition 
as a smooth section of ${\mathbb E}_{\rho}$ 
under the trivialization given by the basis 
$\left\{ du_{1}, ..., du_{g} \right\}$ of 
$H^{1,0} \left( {\cal X}/{\cal H}_{g} \right)$. 
Furthermore, denote by $f'$ the associated ${\mathbb C}^{d}$-valued smooth 
function of $Z$ under the trivialization given by the basis of 
$H^{1,0} \left( {\cal X}/{\cal H}_{g} \right)$ 
which is derived from $\left\{ \beta_{1}, ..., \beta_{g} \right\}$. 
Then $f' = \rho \left( Z - \overline{Z} \right) f$ and 
$\nabla(\beta_{i}) = 0$. 
Since $dz_{ij}$ corresponds to $2 \pi \sqrt{-1} (du_{i} du_{j})$ 
under the Kodaira-Spencer isomorphism, 
\begin{eqnarray*} 
(D_{\rho} f) |_{u_{ij} = \beta_{i} \beta_{j}} 
& = & 
\rho \left( Z - \overline{Z} \right)^{-1} 
C \left( \rho \left( Z - \overline{Z} \right) f \right) 
|_{u_{ij} = \beta_{i} \beta_{j}} 
\\ 
& = & 
\rho \left( Z - \overline{Z} \right)^{-1} \nabla(f') 
\ \ {\rm mod} \left( H^{0,1} \left( {\cal X}/{\cal H}_{g} \right) \right) 
\\ 
& = & 
\nabla(f) 
\ \ {\rm mod} \left( H^{0,1} \left( {\cal X}/{\cal H}_{g} \right) \right), 
\end{eqnarray*} 
and hence the assertion follows from the induction on $e$. 
\ $\square$ 
\vspace{2ex}

Let $\kappa \in X^{+}(T_{g})$ be as above. 
Then the Gauss-Manin connection gives 
$$
{\mathbb D}_{\kappa} \rightarrow 
{\mathbb D}_{\kappa} \otimes \left( \Omega_{{\cal A}_{g,N}}^{1} \right)^{e}. 
$$
This, together with the Kodaira-Spencer isomorphism 
$$
\Omega_{{\cal A}_{g,N}}^{1} \cong 
{\rm Sym}^{2} \left( \pi_{*} \left( \Omega_{{\cal X}/{\cal A}_{g,N}}^{1} \right) \right), 
$$ 
gives rise to 
$$
{\mathbb D}_{\kappa} \rightarrow 
{\mathbb D}_{\kappa} \otimes 
\left( {\rm Sym}^{2} \left( \pi_{*} \left( \Omega_{{\cal X}/{\cal A}_{g,N}}^{1} 
\right) \right) \right)^{e}
$$
which we denote by ${\cal D}_{\kappa}^{e}$. 
\vspace{2ex}  

{\sc Proposition 3.2.} 
\begin{it} 
Let $\rho : GL_{g} \rightarrow GL_{d}$ be the rational homomorphism associated with 
$W_{\kappa}$. 
Then via the projection 
$H_{\rm DR}^{1} \left( {\cal X}/{\cal H}_{g} \right) \rightarrow 
\pi_{*} \left( \Omega_{{\cal X}/{\cal H}_{g}}^{1} \right)$ 
derived from the Hodge decomposition, 
${\cal D}_{\kappa}^{e}$ gives $D_{\rho}^{e}$. 
\end{it} 
\vspace{2ex}

{\it Proof.} 
This assertion follows from Proposition 3.1. \ $\square$ 
\vspace{2ex}

{\bf 3.2. Nearly holomorphic modular forms.} \ 
We recall the definition of nearly holomorphic Siegel modular forms by Shimura. 
\vspace{2ex}

{\sc Definition 3.3.} 
Let $R$ be a ${\mathbb Z} \left[ 1/N, \zeta_{N} \right]$-subalgebra of ${\mathbb C}$. 
A ${\mathbb C}^{d}$-valued smooth function $f$ of 
$Z = X + \sqrt{-1} Y \in {\cal H}_{g}$ is defined to be 
{\it nearly holomorphic over $R$} if $f$ has the following expression 
$$
f(Z) = \sum_{T} q \left( T, \pi^{-1} Y^{-1} \right) \cdot 
\exp \left( 2 \pi \sqrt{-1} {\rm tr}(TZ)/N \right), 
$$
where $q \left( T, \pi^{-1} Y^{-1} \right)$ are vectors of degree $d$ 
whose entries are polynomials over $R$ of the entries of $(4 \pi Y)^{-1}$. 
For a rational homomorphism $\rho : GL_{g} \rightarrow GL_{d}$ over $R$, 
denote by ${\cal N}_{\rho}^{\rm hol}(R)$ the $R$-module of 
all ${\mathbb C}^{d}$-valued smooth functions which are nearly holomorphic over $R$ with $\rho$-automorphic condition for $\Gamma(N)$. 
Call these elements {\it nearly holomorphic Siegel modular forms over $R$ 
of weight $\rho$} (and degree $g$, level $N$). 
\vspace{2ex}

{\sc Theorem 3.4.} 
\begin{it} 
Let $R$ be a ${\mathbb Z} \left[ 1/N, \zeta_{N} \right]$-subalgebra of ${\mathbb C}$, 
and $\rho : GL_{g} \rightarrow GL_{d}$ be a rational homomorphism over $R$ 
associated with $W_{\kappa - h(1,..., 1)}$ for 
$\kappa \in X^{+}(T_{g})$, $h \in {\mathbb Z}$. 
Then there exists a natural $R$-linear isomorphism 
$$
\Phi : {\cal N}_{(\kappa, h)}(R) \rightarrow {\cal N}_{\rho}^{\rm hol}(R). 
$$
Consequently,  ${\cal N}_{\rho}^{\rm hol}(R)$ is a finitely generated $R$-module, 
and 
\end{it} 
$$
{\cal N}_{\rho}^{\rm hol}(R) \otimes_{R} {\mathbb C} = 
{\cal N}_{\rho}^{\rm hol}({\mathbb C}). 
$$

{\it Proof.} 
The last assertion follows from the existence of $\Phi$, 
Theorem 2.3 and that 
${\cal N}_{(\kappa, h)}(R) \otimes_{R} {\mathbb C} = 
{\cal N}_{(\kappa, h)}({\mathbb C})$. 
Therefore, 
we will show the existence of $\Phi$ with required property. 
Since $H^{0,1} \left( {\cal X}/{\cal H}_{g} \right)$ is an anisotropic subspace of 
$H_{\rm DR}^{1} \left( {\cal X}/{\cal H}_{g} \right)$, 
as is seen in 2.1, the projection 
$H_{\rm DR}^{1} \left( {\cal X}/{\cal H}_{g} \right) \rightarrow 
H^{1,0} \left( {\cal X}/{\cal H}_{g} \right)$ 
gives rise to the homomorphism 
${\mathbb D}_{\kappa} \rightarrow {\mathbb E}_{\kappa}$ 
on ${\cal A}_{g,N}({\mathbb C})$, 
and hence one has the ${\mathbb C}$-linear map 
$\Phi : {\cal N}_{(\kappa, h)}({\mathbb C}) \rightarrow 
{\cal N}_{\rho}^{\rm hol}({\mathbb C})$. 

First, we show that 
$\Phi \left( {\cal N}_{(\kappa, h)}(R) \right) \subset {\cal N}_{\rho}^{\rm hol}(R)$ 
and $\Phi$ is injective. 
Since 
$$
d u_{i} - d \overline{u_{i}} = 2 \sqrt{-1} \sum_{j=1}^{g} {\rm Im}(z_{ij}) \beta_{j}, 
$$
if we put 
$W = (w_{ij})_{i,j} = {\rm Im}(Z)^{-1}$, 
then 
$$
\eta_{i} = \beta_{i} = 
\frac{1}{2 \pi \sqrt{-1}} \sum_{j=1}^{g} w_{ij} \left( d u_{i} - d \overline{u_{i}} \right). 
$$ 
Let $X_{ij}$ $(1 \leq i, j \leq 2g)$ be the variables corresponding to $e_{i} \otimes e_{j}$, 
where 
$$
e_{i} = \left\{ \begin{array}{ll} 
\omega_{i} & (1 \leq i \leq g), \\ \eta_{i-g} & (g+1 \leq i \leq 2g). 
\end{array} \right.  
$$
Then under the trivialization of $\det({\mathbb E})^{\otimes h}$ by 
$(\omega_{1} \wedge \cdots \wedge \omega_{g})^{\otimes h}$,  
each $f \in {\cal N}_{(\kappa, h)}({\mathbb C})$ is expressed as a polynomial of 
$X_{ij}$ $(1 \leq i, j \leq 2g)$ over the ring ${\cal O}_{{\cal H}_{g}}^{\rm hol}$ of  holomorphic functions of $Z \in {\cal H}_{g}$. 
Therefore, 
$\Phi(f)$ is a polynomial of $X_{ij}$ $(1 \leq i, j \leq g)$ over the ring of 
smooth functions on ${\cal H}_{g}$ which is obtained from $f$ by substituting 
$$
e_{g+i} = \frac{1}{2 \sqrt{-1}} \sum_{j=1}^{g} 
w_{ij} \left( \frac{\omega_{j}}{2 \pi \sqrt{-1}} - d \overline{u_{j}} \right) = 
\sum_{j=1}^{g} \left( \frac{w_{ij}}{-4 \pi} \omega_{j} - 
\frac{w_{ij}}{2 \sqrt{-1}} d \overline{u_{j}} \right), 
$$
and putting $d \overline{u_{i}} = 0$. 
This process maps $X_{ij}$ $(1 \leq i, j \leq 2g)$ to    
$$
\left\{ \begin{array}{ll} 
X_{ij} & (i, j \leq g), 
\\
{\displaystyle \sum_{k=1}^{g} \frac{w_{i-g, k}}{-4 \pi} X_{kj}} & (i > g, \ j \leq g), 
\\
{\displaystyle \sum_{l=1}^{g} \frac{w_{j-g, l}}{-4 \pi} X_{il}} & (i \leq g, \ j > g), 
\\
{\displaystyle \sum_{k,l=1}^{g} \frac{w_{i-g, k}}{-4 \pi} \frac{w_{j-g, l}}{-4 \pi} X_{kl}} & 
(i > g, \ j > g),  
\end{array} \right. 
$$
where $w_{i,j}$ denotes $w_{ij}$, 
and hence 
$\Phi \left( {\cal N}_{(\kappa, h)}(R) \right) \subset {\cal N}_{\rho}^{\rm hol}(R)$. 
Furthermore, 
if the images of $X_{ij}$ and $X_{kl}$ $(1 \leq i, j, k, l \leq 2g)$ 
have a common nonzero term, 
then $X_{ij} = X_{kl}$. 
Therefore, 
since the functions $w_{ij}$ $(1 \leq i, j \leq g)$ on ${\cal H}_{g}$ are 
algebraically independent over ${\cal O}_{{\cal H}_{g}}^{\rm hol}$, 
$\Phi(f) = 0$ implies that $f = 0$. 
This means the injectivity of $\Phi$. 

Second, we show that $\Phi$ is surjective. 
Let $f$ be an element of ${\cal N}_{\rho}^{\rm hol}(R)$. 
Then by a result of Shimura \cite[14.2. Proposition]{Sh}, 
if a non-zero element $\psi$ of ${\cal M}_{a}({\mathbb C})$ has 
a sufficiently large weight $a$, 
then there are positive integers $n, e_{1},..., e_{n}$ and 
$g_{e_{i}} \in {\cal M}_{\rho \otimes \sigma^{e_{i}}}({\mathbb C})$ such that 
$$
f \cdot \psi = 
\sum_{i=1}^{n} \left( \theta_{e_{i}} \circ 
D_{\rho \otimes \sigma^{e_{i}}}^{e_{i}} \right) \left( g_{e_{i}} \right). 
$$ 
As seen in 2.1, 
${\mathbb E}_{\rho \otimes \sigma^{e_{i}}}$ becomes a direct summand of 
${\mathbb D}_{(\kappa \otimes \sigma^{e_{i}}, h)}$ over ${\mathbb C}$ 
as its holomorphic part, 
and hence $g_{e_{i}}$ can be regarded as a global section of 
${\mathbb D}_{(\kappa \otimes \sigma^{e_{i}}, h)}$ over ${\mathbb C}$. 
This implies that  
$$
\sum_{i=1}^{n} \left( \theta_{e_{i}} \circ 
{\cal D}_{\rho \otimes \sigma^{e_{i}}}^{e_{i}} \right) \left( g_{e_{i}} \right) 
$$ 
is a global section of ${\mathbb D}_{(\kappa, h)} \otimes \det({\mathbb E})^{\otimes a}$ 
over ${\mathbb C}$ whose image by $\Phi$ becomes $f \cdot \psi$. 
Therefore,  
there exists a meromorphic section $f'$ of ${\mathbb D}_{(\kappa, h)}$ 
over ${\cal A}_{g,N}({\mathbb C})$ 
such that $\Phi(f') = f$ and that $f'$ is regular outside the divisor given by $\psi = 0$. 
By the injectivity of $\Phi$, 
$f'$ is independent of the choice of $\psi$. 
We will show that $f'$ is regular on ${\cal A}_{g,N}$. 
Since $\omega^{*}$ is ample on ${\cal A}^{*}_{g,N}$, 
there are positive integers $c$ (: sufficiently large) and $d$ such that 
$(\omega^{*})^{\otimes c}$ gives rise to a closed immersion 
$\iota : {\cal A}^{*}_{g,N} \hookrightarrow {\mathbb P}^{d}$ 
and that 
$\iota^{*} \left( {\cal O}_{{\mathbb P}^{d}}(1) \right) \cong (\omega^{*})^{\otimes c}$.   
Under the representation of points on ${\mathbb P}^{d}$ 
by homogeneous coordinates $[x_{0}: \cdots : x_{d}]$, 
one has $\iota^{*}(x_{i}) \in {\cal M}_{c}({\mathbb C})$. 
Then by putting $\psi = \iota^{*}(x_{i})$, 
one can see that $f'$ is regular where $\iota^{*}(x_{i}) \neq 0$, 
and hence  $f'$ is regular on the whole ${\cal A}_{g,N}$. 
Since the Fourier coefficients of $\Phi(f') = f \in {\cal N}_{\rho}^{\rm hol}(R)$ are 
polynomials over $R$ of $X_{ij}$ and $w_{ij}/(4 \pi)$ $(1 \leq i, j \leq g)$, 
those of $f' \in {\cal N}_{(\kappa, h)}({\mathbb C})$ 
are polynomials over $R$ of  $X_{ij}$ $(1 \leq i, j \leq 2g)$, 
and hence by Theorem 2.4, 
one has $f' \in {\cal N}_{(\kappa, h)}(R)$. \ $\square$ 
\vspace{2ex}

{\sc Theorem 3.5.} 
\begin{it} 
Let $R$ be a ${\mathbb Z} \left[ 1/N, \zeta_{N} \right]$-subalgebra of ${\mathbb C}$, 
and $\rho : GL_{g} \rightarrow GL_{d}$ be a rational homomorphism over $R$ 
associated with $W_{\kappa - h(1,..., 1)}$ as in Theorem 3.4. 
Let $\widetilde{\alpha}$ be a test object over $R$ corresponding to 
a CM abelian variety $X$, 
and assume that one can extend the basis of $\Omega_{X/R}^{1}$ to 
a basis of $H^{1}(X/R)$ which gives a projection 
$H^{1}(X/R) \rightarrow \Omega_{X/R}^{1}$ compatible with the action of 
${\rm End}(X)$. 
Then for any $f \in {\cal N}_{\rho}^{\rm hol}(R)$, 
the evaluation $f \left( \widetilde{\alpha} \right)$ of $f$ 
at $\widetilde{\alpha}$ belongs to $R^{d}$. 
\end{it}
\vspace{2ex}

{\it Proof.} \ 
This assertion follows from Theorems 2.5 and 3.4. \ $\square$ 
\vspace{3ex}

{\bf 4. Arithmeticity in the {\boldmath $p$}-adic case} 
\vspace{2ex}

{\bf 4.1. {\boldmath $p$}-adic modular forms.} \ 
According to \cite{S}, 
we consider $p$-adic modular forms as limits of modular forms 
defined in 2.4. 
In what follows, fix a prime $p$ not dividing $N \geq 3$, 
and a $p$-adic field, 
i.e., a finite extension $K$ of ${\mathbb Q}_{p}$. 
Furthermore, assume that $K$ contains $\zeta_{N}$. 
Then the valuation ring $R_{K}$ of $K$ contains $1/N$ and $\zeta_{N}$. 
Fix a generator $\pi$ of the maximal ideal of $R_{K}$, 
and put $R_{m} = R_{K} / (\pi^{m})$ for positive integers $m$. 

{\it Weights} of vector-valued $p$-adic Siegel modular forms over $R_{K}$ 
are defined as continuous homomorphisms 
$\rho : GL_{g}({\mathbb Z}_{p}) \rightarrow GL_{d}(R_{K})$ 
such that there are rational homomorphisms 
$\rho_{k} : GL_{g} \rightarrow GL_{d}$ over $R_{K}$ 
with $\rho = \lim_{k} \rho_{k}$ which means that $\rho_{k}$ converge to $\rho$ 
uniformly on $GL_{g}({\mathbb Z}_{p})$. 
Then $\rho(\alpha) \ {\rm mod}(\pi^{m})$ 
$\left( \alpha \in GL_{g}({\mathbb Z}_{p}) \right)$ 
is given by a rational homomorphism $GL_{g} \rightarrow GL_{d}$ over $R_{m}$ 
which we denote by $\rho(m)$. 
Denote by ${\cal W}_{d}(R_{K})$ the set of these weights. 

For $\rho = \lim_{k} \rho_{k} \in {\cal W}_{d}(R_{K})$, 
a sequence $\{ f_{k} \}_{k}$ of $f_{k} \in {\cal M}_{\rho_{k}}(R_{K})$ is 
a Cauchy sequence if the following condition holds: 
for each $m$, 
there is a positive integer $k(m)$ such that if $k \geq k(m)$, 
then $\rho \equiv \rho_{k}$ ${\rm mod}(\pi^{m})$ and the image of $f_{k}$ by 
the reduction map 
${\cal M}_{\rho}(R_{K}) \rightarrow {\cal M}_{\rho(m)}(R_{m})$ 
is independent of $k$. 
Two Cauchy sequences 
$\left\{ f_{k} \in {\cal M}_{\rho_{k}}(R_{K}) \right\}_{k}$ and 
$\left\{ f'_{k} \in {\cal M}_{\rho'_{k}}(R_{K}) \right\}_{k}$  
with $\rho = \lim_{k} \rho_{k} = \lim_{k} \rho'_{k}$ are called equivalent 
if the images of $f_{k}$ and $f'_{k}$ in 
${\cal M}_{\rho(m)}(R_{m})$ have the same limit 
for any positive integer $m$. 
\vspace{2ex}

{\sc Definition 4.1.} 
Let $K$ be a $p$-adic field, 
and $\rho = \lim_{k} \rho_{k} \in {\cal W}_{d}(R_{K})$ 
be a weight over the valuation ring $R_{K}$ of $K$. 
Then the $R_{K}$-module 
$\overline{\cal M}_{\rho}(R_{K})$ 
of {\it $p$-adic Siegel modular forms over $R_{K}$ of weight $\rho$} 
(and degree $g$, level $N$) is defined as the equivalence classes of 
the above Cauchy sequences. 
We also put 
$$
\overline{\cal M}_{\rho}(K) = 
\overline{\cal M}_{\rho}(R_{K}) \otimes K 
$$
whose elements 
$f = \lim_{k} f_{k}$ $\left( f_{k} \in {\cal M}_{\rho_{k}}(K) \right)$ 
are called {\it $p$-adic Siegel modular forms over $K$ of weight $\rho$} 
(and degree $g$, level $N$). 
\vspace{2ex}

Let $c$ be a $0$-dimensional cusp on ${\cal A}_{g,N}^{\ast}$. 
Then there are trivializations over 
${\rm Spec} \left( {\cal R}_{g,N} \otimes R_{m} \right) = 
{\rm Spec} \left( {\cal R}_{g,N} 
\otimes_{{\mathbb Z}[1/N, \zeta_{N}]} R_{m} \right)$: 
$$
{\mathbb E}_{\rho(m)} \times_{{\cal A}_{g,N} \otimes R_{m}} 
{\rm Spec} \left( {\cal R}_{g,N} \otimes R_{m} \right) = 
\left( {\cal R}_{g,N} \otimes R_{m} \right)^{d} 
$$ 
compatible with $m$, 
from which one has the Fourier expansion map 
$$
F_{c} : \overline{\cal M}_{\rho}(K) \rightarrow 
\left( {\cal R}_{g,N} \otimes_{{\mathbb Z}[1/N, \zeta_{N}]} K \right)^{d} 
$$
as $F_{c}(f) = \lim_{k} F_{c} \left( f_{k} \right)$. 
By the $q$-expansion principle, 
${\cal M}_{\rho_{k}}(R_{K})$  consists of elements of 
${\cal M}_{\rho_{k}}(K)$ with Fourier coefficients in $R_{K}$, 
and hence 
$$
\overline{\cal M}_{\rho}(R_{K}) = 
\left\{ f \in \overline{\cal M}_{\rho}(K) 
\ \left| \ F_{c}(f) \in 
\left( {\cal R}_{g,N} \otimes_{{\mathbb Z}[1/N, \zeta_{N}]} R_{K}  \right)^{d} 
\right. \right\}. 
$$
It is clear that if $\rho$ is a rational homomorphism, 
then there are natural inclusions 
${\cal M}_{\rho}(R_{K}) \hookrightarrow 
\overline{\cal M}_{\rho}(R_{K})$ and 
${\cal M}_{\rho}(K) \hookrightarrow \overline{\cal M}_{\rho}(K)$. 
\vspace{2ex} 

{\bf 4.2. Igusa tower.} \ 
We give a description of vector-valued $p$-adic Siegel modular forms 
as automorphic functions on the Igusa tower. 
A similar consideration is given by Hida \cite{Hi}. 
Let $U_{m}$ be the ordinary locus of ${\cal A}_{g,N} \otimes R_{m}$ 
which is a nonempty open subset, and hence is irreducible. 
Denote by ${\cal X}$ the universal abelian scheme over $U_{m}$. 
Then the maximal \'{e}tale quotient ${\cal X}[p^{n}]^{\rm et}$ of 
$$
{\cal X}[p^{n}] = {\rm Ker} 
\left( p^{n} : {\cal X} \rightarrow {\cal X} \right) 
$$
is an \'{e}tale sheaf on $U_{m}$ 
of free $\left( {\mathbb Z}/(p^{n}) \right)$-modules of rank $g$, 
and hence one has the associated monodromy representation 
$$
\mu_{m,n} : \pi_{1} \left( U_{m} \right) \rightarrow 
GL_{g} \left( {\mathbb Z}/(p^{n}) \right). 
$$ 
Then by a result of Chai and Faltings \cite[Chapter V, 7.2]{FC}, 
$\mu_{m,n}$ is surjective, 
and hence by the irreducibility of $U_{m}$, 
the Galois covering $T_{m,n}$ of $U_{m}$ associated with 
${\rm Ker} \left( \mu_{m,n} \right)$ is irreducible. 
The system $\{ T_{m,n} \}_{m,n}$ is called the {\it Igusa tower} 
which satisfies that 
$$
T_{m+1,n} \otimes R_{m} \cong T_{m,n}. 
$$
In particular, 
$T_{m+1,m+1} \otimes R_{m}$ is a Galois covering of $T_{m,m}$ 
with Galois group 
$$
{\rm Ker} \left( GL_{g} \left( {\mathbb Z}/(p^{m+1}) \right) \rightarrow 
GL_{g} \left( {\mathbb Z}/(p^{m}) \right) \right). 
$$
Then an $R_{K}^{d}$-valued function on $\{ T_{m,m} \}_{m}$ is defined as 
a set of elements 
$$ 
\phi_{m} \in H^{0} 
\left( T_{m,m}, {\cal O}_{T_{m,m}} \otimes (R_{m})^{d} \right) 
$$
such that the restriction of $\phi_{m+1}$ to $T_{m+1,m+1} \otimes R_{m}$ 
is reduced to $\phi_{m}$ under the projection 
$T_{m+1,m+1} \otimes R_{m} \rightarrow T_{m,m} \otimes R_{m}$. 
Furthermore, the $\rho$-automorphic condition means that 
$$
\rho(m)(\alpha) \cdot \alpha(\phi_{m}) = \phi_{m} 
$$
for any $\alpha \in GL_{g}({\mathbb Z}_{p})$. 
As seen in 2.5, for a $0$-dimensional cusp $c$ on ${\cal A}^{*}_{g,N}$, 
the Mumford uniformization theory gives 
the associated principally polarized abelian scheme 
${\cal X}_{c}$ with symplectic level $N$ structure over ${\cal R}_{g,N}$. 
Since any geometric fiber of 
${\cal X}_{c} \otimes_{{\mathbb Z}[1/N, \zeta_{N}]} R_{m}$ is ordinary, 
${\cal X}_{c}$ gives rise to a principally polarized abelian scheme 
with symplectic level $N$ structure over $\{ T_{m,m} \}$, 
and the Fourier expansion map $F_{c}$ on $\overline{\cal M}_{\rho}(R)$ 
becomes the evaluation on this abelian scheme. 
Then the following characterization of $p$-adic Siegel modular forms 
is given in \cite[Theorem 4]{I}: 
\vspace{2ex}

\begin{it} 
The $R_{K}$-module $\overline{\cal M}_{\rho}(R_{K})$ of 
$p$-adic Siegel modular forms over $R_{K}$ of weight $\rho$ 
consists of $R_{K}^{d}$-valued functions $\phi = \{ \phi_{m} \}$ 
on the Igusa tower $\{ T_{m,m} \}$ with $\rho$-automorphic condition that 
$$
F_{c} \left( \alpha(\phi) \right) = \rho(\alpha)^{-1} \cdot F_{c}(\phi), 
$$
where $\alpha \in GL_{g}({\mathbb Z}_{p})$ acts on 
$\overline{\cal M}_{\rho}(R_{K})$ 
via its action on $\{ T_{m,m} \}$. 
\end{it}
\vspace{2ex}

{\bf 4.3. {\boldmath $p$}-adic differential operator.} \ 
We give a $p$-adic counterpart of Shimura's differential operator. 
\vspace{2ex}

{\sc Proposition 4.2.} 
\begin{it} 
Denote by $R_{K}$ the valuation ring of a $p$-adic field $K$ 
containing $\zeta_{N}$. 

{\rm (1)} 
For each weight $\rho \in {\cal W}_{d}(R_{K})$, 
there exists an $R_{K}$-linear map 
$$
D_{p, \rho}^{e} : \overline{\cal M}_{\rho}(R_{K})  \rightarrow 
\overline{\cal M}_{\rho \otimes \tau^{e}}(R_{K})
$$
which is defined inductively as 
$$
F_{c} \left( D_{p, \rho}^{e}(f) \right) = \sum_{1 \leq i \leq j \leq g} q_{ij} 
\frac{\partial F_{c} \left( D_{p, \rho}^{e-1}(f) \right)}{\partial q_{ij}} 
$$ 
for each $f \in \overline{\cal M}_{\rho}(R_{K})$. 
We extend $D_{p, \rho}^{e}$ to a $K$-linear map 
$$
\overline{\cal M}_{\rho}(K) \rightarrow 
\overline{\cal M}_{\rho \otimes \tau^{e}}(K)
$$ 
which we denote by the same symbol. 

{\rm (2)} 
Let ${\cal U}$ be the formal scheme given by the inverse limit of $U_{m}$, 
and let $\pi : {\cal X} \rightarrow {\cal U}$ be the canonical family of 
abelian schemes. 
Then for each rational homomorphism 
$\rho : GL_{g} \rightarrow GL_{d}$ over $R_{K}$, 
$D_{p, \rho}^{e}$ is obtained from the composition 
$$
{\mathbb E}_{\rho} 
\rightarrow 
{\mathbb E}_{\rho} \otimes 
\left( \Omega_{\cal U}^{1} \right)^{\otimes e} 
\rightarrow 
{\mathbb E}_{\rho} \otimes 
\left( {\rm Sym}^{2} \left( \pi_{*} \left( 
\Omega_{{\cal X}/{\cal U}}^{1} \right) \right) \right)^{\otimes e}. 
$$ 
Here the first map is given by the Gauss-Manin connection 
$$
\nabla : H_{\rm DR}^{1} \left( {\cal X}/{\cal U} \right) \rightarrow 
H_{\rm DR}^{1} \left( {\cal X}/{\cal U} \right) \otimes \Omega_{\cal U}^{1} 
$$
together with the projection 
$H_{\rm DR}^{1} \left( {\cal X}/{\cal U} \right) \rightarrow 
\pi_{*} \left( \Omega^{1}_{{\cal X}/{\cal U}} \right)$ 
derived from the unit root splitting of 
$H_{\rm DR}^{1} \left( {\cal X}/{\cal U} \right)$ (cf. \cite[1.11]{Kt3}), 
and the second map is given by the Kodaira-Spencer isomorphism 
$\Omega_{\cal U}^{1} \cong {\rm Sym}^{2} 
\left( \pi_{*} \left( \Omega_{{\cal X}/{\cal U}}^{1} \right) \right)$. 
\end{it} 
\vspace{2ex}

{\it Proof.} \ 
First, we prove (1). 
Recall that ${\cal X}_{c}$ denotes the principally polarized abelian scheme 
with symplectic level $N$ structure over ${\cal R}_{g,N}$ 
which is associated with a $0$-dimensional cusp $c$. 
Then it is known (cf. \cite[Chapter III, 9]{FC}) that 
the Kodaira-Spencer isomorphism gives 
$$
\frac{d q_{ij}}{q_{ij}} \leftrightarrow \omega_{i} \omega_{j} = 
\frac{d x_{i}}{x_{i}} \frac{d x_{j}}{x_{j}} \ (1 \leq i, j \leq g) 
$$
for ${\cal X}_{c}$, 
and hence $D_{p, \rho}^{e}(f)$ is an 
$S_{e} \left( {\rm Sym}^{2} \left( R_{K}^{g} \right), 
R_{K}^{d} \right)$-valued 
function on the Igusa tower with $\rho$-automorphic condition. 
Therefore, by the characterization of $p$-adic modular forms given in 4.2, 
$D_{p, \rho}^{e}(f) \in 
\overline{\cal M}_{\rho \otimes \tau^{e}}(R_{K})$. 

Second, we prove (2). 
Let 
$$
H_{\rm DR}^{1} \left( {\cal X}/{\cal U} \right) = 
\pi_{*} \left( \Omega_{{\cal X}/{\cal U}}^{1} \right) \oplus 
U_{{\cal X}/{\cal U}} 
$$
be the unit root splitting, 
where $U_{{\cal X}/{\cal U}}$ denotes the unit root subspace for the Frobenius action. 
Then as is shown in \cite[(1.12.7) Key Lemma]{Kt3} and \cite[Lemma 5.9]{E1}, 
from the description of $U_{{\cal X}/{\cal U}}$ 
for the abelian scheme over ${\cal U}$ given by ${\cal X}_{c}$, 
the derivation $D_{ij} = q_{ij} \partial/\partial q_{ij}$ satisfies that 
$$
\nabla(D_{ij}) 
\left( \pi_{*} \left( \Omega_{{\cal X}/{\cal U}}^{1} \right) \right) 
\subset U_{{\cal X}/{\cal U}}. 
$$
Therefore, 
identifying $f \in \overline{\cal M}_{\rho}(R_{K})$ and 
its Fourier expansion, 
one has  
\begin{eqnarray*}
\nabla(f) 
& = & 
\sum_{1 \leq i \leq j \leq g} \nabla \left( D_{ij}(f) \right) 
\left( \frac{dx_{i}}{x_{i}} \right) \left( \frac{dx_{j}}{x_{j}} \right) 
\\
& = & 
\sum_{1 \leq i \leq j \leq g} q_{ij} \frac{\partial f}{\partial q_{ij}} 
\left( \frac{dx_{i}}{x_{i}} \right) \left( \frac{dx_{j}}{x_{j}} \right) 
\ {\rm mod} \left( U_{{\cal X}/{\cal U}} \right), 
\end{eqnarray*}
and hence 
$$
D_{p, \rho}(f) = 
\sum_{1 \leq i \leq j \leq g} q_{ij} \frac{\partial f}{\partial q_{ij}}. 
$$
This implies (2) by the induction on $e$. 
\ $\square$ 
\vspace{2ex}

{\it Remark 2.} 
Let 
$$
\phi : K^{d} \cong {\rm Hom}_{K} \left( K, K^{d} \right) \rightarrow 
S_{g} \left( {\rm Sym}^{2} \left( K^{g} \right), K^{d} \right) 
$$
be the $K$-linear map induced from the determinant map 
${\rm Sym}^{2} \left( K^{g} \right) \rightarrow K$. 
Then $\phi$ is $GL_{g}$-equivariant for the representations 
$\rho \otimes \det^{\otimes 2}$ and $\rho \otimes \tau^{g}$, 
and the pullback by $\phi$ gives rise to a $p$-adic differential operator 
$$
\theta: \overline{\cal M}_{\rho}(K) \rightarrow 
\overline{\cal M}_{\rho \otimes \det^{\otimes 2}}(K). 
$$
This is called the theta operator which sends $f$ with Fourier expansion 
$$\sum_{T} a(T) \mbox{\boldmath $q$}^{T/N}
$$
to $\theta(f)$ with Fourier expansion 
$\sum_{T} a(T) \det(T/N) \mbox{\boldmath $q$}^{T/N}$. 
This operator was studied by B\"{o}cherer and Nagaoka \cite{BN1,BN2,BN3}. 
\vspace{2ex}

{\bf 4.4. Correspondence between nearly and {\boldmath $p$}-adic modular forms.} \ 
Let $k$ be a subfield of ${\mathbb C}$ containing $\zeta_{N}$. 
A CM point over $k$ associated with an abelian variety $X$ is called 
a {\it $p$-ordinary CM point} if $X$ has good and ordinary reduction 
at a fixed place of $k$ lying above $p$. 
\vspace{2ex}

{\sc Theorem 4.3.} 
\begin{it} 
Let $K$ be a $p$-adic field containing $k$, 
and $\rho : GL_{g} \rightarrow GL_{d}$ be a rational homomorphism over $k \cap R_{K}$. 
Then there exists uniquely an injective $k$-linear map 
$$
\iota_{p} : {\cal N}_{\rho}^{\rm hol}(k) \rightarrow \overline{\cal M}_{\rho}(K) 
$$
such that for any $f \in {\cal N}_{\rho}^{\rm hol}(k)$ and 
any extended test object $\widetilde{\alpha}$ over $k' \supset k$ 
associated with a $p$-ordinary CM point $\alpha$ 
on ${\cal A}_{g,N}(k')$, 
$$
f \left( \widetilde{\alpha} \right) = 
\iota_{p}(f) \left( \widetilde{\alpha} \right) 
$$
as an element of $(k')^{d}$. 
\end{it}
\vspace{2ex}

{\it Proof.} \ 
We prove the assertion by constructing $\iota_{p}$ and showing its properties. 

{\bf Construction of {\boldmath $\iota_{p}$}:} 
Let $h_{p-1}$ denote the generalized Hasse invariant which is 
the unique Siegel modular form over ${\mathbb F}_{p}$ of degree $g$, 
weight $p-1$ and level $1$ whose Fourier expansion is the constant $1$. 
Then by the ampleness of $\omega^{*}$ on ${\cal A}_{g,N}^{*}$, 
there is a positive integer $a$ and 
$$
\psi \in {\cal M}_{(p-1)a} \left( k \cap R_{K} \right) = 
H^{0} \left( {\cal A}_{g,N}^{*} \otimes \left( k \cap R_{K} \right), 
\left( \omega^{*} \right)^{(p-1)a} \right) 
$$ 
whose reduction is $\left( h_{p-1} \right)^{a}$. 
Let $f$ be an element of ${\cal N}_{\rho}^{\rm hol}(k)$. 
Then by \cite[14.2. Proposition]{Sh}, 
taking a sufficiently large positive integer $b$, 
there are positive integers $n, e_{1},..., e_{n}$ and 
$g_{e_{i}} \in {\cal M}_{\rho \otimes \sigma^{e_{i}}}(k)$ such that
$$
f \cdot \psi^{b} = \sum_{i=1}^{n} 
\left( \theta_{e_{i}} \circ {D}_{\rho \otimes \sigma^{e_{i}}}^{e_{i}} \right) 
\left( g_{e_{i}} \right). 
$$
Then using Proposition 4.2, we put 
$$
\iota_{p}(f) = \psi^{-b} \sum_{i=1}^{n} 
\left( \theta_{e_{i}} \circ 
{D}_{p, \rho \otimes \sigma^{e_{i}}}^{e_{i}} \right) \left( g_{e_{i}} \right) 
$$
which is well defined since $\psi$ is invertible as a $p$-adic modular form 
over $R_{K}$. 

{\bf Preserving {\boldmath $p$}-ordinary CM values 
for {\boldmath $\iota_{p}$}:} 
Consider $H_{\rm DR}^{1} (X/k')$ for the abelian variety 
$X$ over $k'$ corresponding to $\alpha$. 
Then it is shown that in \cite[(5.1.27) Key Lemma]{Kt3} that 
the decomposition  
$$
H_{\rm DR}^{1}(X/k') = 
H_{\rm DR}^{1}(X/k')^{+} \oplus H_{\rm DR}^{1}(X/k')^{-}; \ 
H_{\rm DR}^{1}(X/k')^{+} = \Omega_{X/k'}^{1},  
$$
which is stable under the action of ${\rm End}(X)$, 
gives rise to the Hodge decomposition and unit root splitting of 
$H_{\rm DR}^{1}(X/k') \otimes (K k')$. 
Therefore, by the construction of $\iota_{p}$, 
the evaluation of $f$ and $\iota_{p}(f)$ at $\widetilde{\alpha}$ 
are equal. 

{\bf Well-definedness and uniqueness of {\boldmath $\iota_{p}$}:} 
This is derived from the following fact: 
Let $R$ be a complete discrete valuation ring whose residue field 
$R / {\mathfrak m}_{R}$ is an algebraically closed field 
of characteristic $p$, 
and $X_{0}$ be a principally polarized ordinary abelian variety 
of dimension $g$ over $R / {\mathfrak m}_{R}$. 
Then by Serre-Tate's local moduli theory 
(cf. \cite[Appendix]{Me} and \cite[8.2]{Hi}), 
the liftings of $X_{0}$ as principally polarized abelian schemes over $R$ 
are parametrized by the multiplicative group 
$$
G(R) = 
{\rm Hom}_{{\mathbb Z}_{p}} \left( {\rm Sym}^{2} \left( T_{p}(X_{0}) \right), 
\widehat{\mathbb G}_{m}(R) \right) 
\cong 
\left( 1 + {\mathfrak m}_{R} \right)^{g(g+1)/2}, 
$$
where $T_{p}(X_{0}) \cong {\mathbb Z}_{p}^{\oplus g}$ 
denotes the $p$-adic Tate module of $X_{0}$. 
Furthermore, by the assumption that $p$ is prime to $N$, 
any symplectic level $N$ structure on $X_{0}$ can be lifted uniquely to 
that on its liftings. 
In this parametrization, 
the lifting of $X_{0}$ corresponding to the identity element 
is called the canonical lifting $X^{\rm can}$, 
and the liftings over $R \left[ \zeta_{p^{n}} \right]$ corresponding to 
torsion points on $G$ are called the quasi-canonical liftings. 
Then as is shown in \cite[Appendix]{Me}, 
the canonical lifting is the unique lifting of $X_{0}$ such that 
all endomorphisms of $X_{0}$ lift to $X^{\rm can}$, 
and the quasi-canonical liftings of $X_{0}$ are mutually isogeneous, 
hence of CM type. 
Since $1 \in G$ is a limit point of the set of torsion points 
on $\bigcup_{n} G \left( R[\zeta_{p^{n}}] \right)$, 
by the $p$-adic Weierstrass preparation theorem, 
a $p$-adic function on $\bigcup_{n} G \left( R[\zeta_{p^{n}}] \right)$ 
vanishing at all ordinary CM points becomes $0$. 

{\bf Injectivity of {\boldmath $\iota_{p}$}:} 
This follows from that the Hecke orbit of a point is dense in 
${\cal A}_{g,N}({\mathbb C})$ in the usual topology. 
\ $\square$ 
\vspace{2ex}

Let $R$ be a ${\mathbb Z} \left[ 1/N, \zeta_{N} \right]$-subalgebra of 
a $p$-adic field $K$, 
and $\rho : GL_{g} \rightarrow GL_{d}$ be the rational homomorphism over $R$ 
associated with $W_{\kappa - h(1,..., 1)}$ 
for $\kappa \in X^{+}(T_{g})$, $h \in {\mathbb Z}$. 
Since the canonical perfect pairing 
$$
H_{\rm DR}^{1} \left( {\cal X}/{\cal U} \right) \times 
H_{\rm DR}^{1} \left( {\cal X}/{\cal U} \right) \rightarrow 
H_{\rm DR}^{2} \left( {\cal X}/{\cal U} \right)
$$
is equivariant for the Frobenius action, 
the unit root subspace $U_{{\cal X}/{\cal U}}$ is anisotropic in 
$H_{\rm DR}^{1} \left( {\cal X}/{\cal U} \right)$. 
Hence by the associated projection 
$H_{\rm DR}^{1} \left( {\cal X}/{\cal U} \right) \rightarrow 
\pi_{*} \left( \Omega_{{\cal X}/{\cal U}}^{1} \right)$, 
one has an $R$-linear map 
$$
\Phi_{p} : {\cal N}_{(\kappa, h)}(R) \rightarrow \overline{\cal M}_{\rho}(K). 
$$
 
{\sc Theorem 4.4.} 
\begin{it} 
Assume that $R \subset {\mathbb C}$. 
Then $\Phi_{p}$ satisfies $\Phi_{p} = \iota_{p} \circ \Phi$, 
and it is injective. 
\end{it} 
\vspace{2ex}

{\it Proof.} \ 
This assertion follows from the construction and injectivity of $\iota_{p}$ shown in Theorem 4.3, 
and the injectivity of $\Phi$ shown in Theorem 3.4. \ $\square$ 
\vspace{2ex}

{\bf 4.5 Siegel-Eisenstein series.} \ 
Let $\chi$ be a Dirichlet character modulo a positive integer $M$, 
and define the Siegel-Eisenstein series as a function of 
$Z = X + \sqrt{-1} Y \in {\cal H}_{g}$ as 
$$
E_{h}(Z, s, \chi) = 
\det(Y)^{s} \sum_{\gamma \in (P \cap \Gamma_{0}(M)) \backslash \Gamma_{0}(M)} 
\chi \left( \det(D_{\gamma}) \right) j(\gamma, Z)^{-h} 
\left| j(\gamma, Z) \right|^{-2s}, 
$$
where 
\begin{eqnarray*} 
\Gamma_{0}(M) 
& = & 
\left\{ \left. \gamma = \left( \begin{array}{cc} 
A_{\gamma} & B_{\gamma} \\ C_{\gamma} & D_{\gamma} 
\end{array} \right) \in Sp_{2g}({\mathbb Z}) \ \right| \ 
C_{\gamma} \equiv 0 \ {\rm mod}(M) \right\}, 
\\
P 
& = & 
\left\{ \left( \begin{array}{cc} A & B \\ 0 & D \end{array} \right) 
\in Sp_{2g}({\mathbb R}) \right\}. 
\end{eqnarray*}
Then $E(Z, s)$ is absolutely convergent for ${\rm Re}(s) > (g+1-h)/2$ 
and analytically continued to the whole complex $s$ plane. 
It is known (cf. \cite[\S 19]{M}) that 
$$
\delta_{h} \left( E_{h}(Z, s, \chi) \right) 
= \pi^{-g} \varepsilon_{g}(h) E_{h+2}(Z, s-1, \chi), 
$$
where $\delta_{h}$ is the Maass-Shimura differential operator given by 
$$
\det \left( Z - \overline{Z} \right)^{((g-1)/2) - h} 
\det \left( 
\frac{1 + \delta_{ij}}{2} \frac{\partial}{\partial (2 \pi \sqrt{-1} z_{ij})} 
\right)_{i,j} 
\det \left( Z - \overline{Z} \right)^{h - ((g-1)/2)}, 
$$
and 
$$
\varepsilon_{g}(h) = 
h \left( h - \frac{1}{2} \right) \cdots \left( h - \frac{g-1}{2} \right). 
$$ 
Furthermore, it is shown in \cite{Si} and \cite{K} that if $h > g+1$, 
then $E_{h}(Z, 0, \chi)$ is a (holomorphic) Siegel modular form 
of degree $g$, weight $h$ and level $M$ whose Fourier coefficients belong to 
a finite cyclotomic extension $k$ of ${\mathbb Q}$. 
Therefore, if $s$ is a nonpositive integer such that $h + 2s > g + 1$, 
then 
$$
\pi^{gs} E_{h}(Z, s, \chi) = 
\prod_{i=0}^{-s-1} \varepsilon_{g}(h + 2s + 2i)^{-1} 
\left( \delta_{h-2} \circ \delta_{h-4} \circ \cdots \circ \delta_{h+2s} 
\right) \left( E_{h+2s}(Z, 0, \chi) \right)
$$
is a nearly holomorphic Siegel modular form of degree $g$, weight $h$ 
and level $M$ which is defined over $k$. 
This Fourier expansion is calculated by Feit, Katsurada, 
Mizumoto and Shimura \cite{F,Ka,Mi,Sh}. 
\vspace{2ex}

{\sc Theorem 4.5.} 
\begin{it}
Let $N \geq 3$ be a multiple of $M$, 
$p$ be a prime not dividing $N$. 
Then for integers $h, s$ satisfying that $(g + 1 - h)/2 < s \leq 0$, 
$$
\iota_{p} \left( \pi^{gs} E_{h}(Z, s, \chi) \right) 
= \prod_{i=0}^{-s-1} \varepsilon_{g}(h + 2s + 2i)^{-1} 
\sum_{T} b_{h + 2s}(T) \det(T)^{-s} \mbox{\boldmath $q$}^{T}, 
$$
where 
$$
E_{h+2s}(Z, 0, \chi) = \sum_{T} b_{h + 2s}(T) \mbox{\boldmath $q$}^{T}. 
$$
\end{it}

{\it Proof.} \ 
Since 
$$
\delta_{h} = 
\left( {\rm id}_{\det({\mathbb E})^{\otimes h}} \otimes \det \right) \circ 
D_{\det^{\otimes h}}, 
$$
its $p$-adic counterpart  
$\left( {\rm id}_{\det({\mathbb E})^{\otimes h}} \otimes \det \right) \circ 
D_{p, \det^{\otimes h}}$ 
is the above theta operator $\theta$ which sends 
$\sum_{T} a(T) \mbox{\boldmath $q$}^{T}$ to 
$\sum_{T} a(T) \det(T) \mbox{\boldmath $q$}^{T}$. 
Then by the construction of $\iota_{p}$ given in the proof of Theorem 4.3, 
\begin{eqnarray*} 
\iota_{p} \left( \pi^{gs} E_{h}(Z, s, \chi) \right) 
& = & 
\prod_{i=0}^{-s-1} \varepsilon_{g}(h + 2s + 2i)^{-1} \cdot 
\theta^{-s} \left( E_{h+2s}(Z, 0, \chi) \right) 
\\ 
& = & 
\prod_{i=0}^{-s-1} \varepsilon_{g}(h + 2s + 2i)^{-1} 
\sum_{T} b_{h + 2s}(T) \det(T)^{-s} \mbox{\boldmath $q$}^{T}. 
\end{eqnarray*} 
This completes the proof. 
\ $\square$ 
\vspace{2ex}

{\bf 4.6. Nearly overconvergence.} \ 
Let $({\cal A}_{g,N})_{\rm rig}$ be the $p$-adic rigid analytic space over $R_{K}$ 
associated with  ${\cal A}_{g,N} \otimes R_{K}$ 
which contains the subspace $({\cal A}_{g,N})_{\rm ord}$ 
associated with the ordinary locus. 
As in the proof of Theorem 4.3, 
there are a positive number $a$ and a lift $\psi \in {\cal M}_{(p-1)a}(R_{K})$ 
of the $a$th power of the generalized Hasse invariant $h_{p-1}$. 
Then for $t \in p^{\mathbb Q} \cap \left[ p^{-1/(p+1)}, 1 \right]$, 
let $({\cal A}_{g,N})_{\rm rig}^{\geq t}$ be the rigid subspace of $({\cal A}_{g,N})_{\rm rig}$ 
defined as the set $x \in ({\cal A}_{g,N})_{\rm rig}$ satisfying $|\psi(x)|_{p} \geq t^{a}$, 
where $|p|_{p} = 1/p$. 
Then for $\kappa \in X^{+}(T_{g})$, $h \in {\mathbb Z}$, there exist natural maps 
$$
H^{0} \left( ({\cal A}_{g,N})_{\rm rig}^{\geq t}, {\mathbb D}_{\kappa} \right) \rightarrow 
H^{0} \left( ({\cal A}_{g,N})_{\rm ord}, {\mathbb D}_{\kappa} \right) \rightarrow 
\overline{\cal M}_{\rho}(K), 
$$
where $\rho$ is associated with $W_{\kappa - h(1,..., 1)}$, 
and we denote this composite by $\varphi_{t}$. 
We define 
$$
{\cal N}_{(\kappa, h)}^{\dagger}(R_{K}) = \bigcup_{t < 1} {\rm Im}(\varphi_{t}), 
$$ 
and call this elements 
{\it nearly overconvergent $p$-adic Siegel modular forms over $R_{K}$ 
of weight $(\kappa, h)$} (of degree $g$, level $N$). 
Since 
${\cal N}_{(\kappa, h)}(K) = {\cal N}_{(\kappa, h)}(R_{K}) \otimes_{R_{K}} K$, 
$\Phi_{p}$ gives rise to a $K$-linear map from ${\cal N}_{(\kappa, h)}(K)$ into the space 
$$
{\cal N}_{(\kappa, h)}^{\dagger}(K) = 
{\cal N}_{(\kappa, h)}^{\dagger}(R_{K}) \otimes_{R_{K}} K 
$$
of nearly overconvergent $p$-adic Siegel modular forms over $K$ 
of weight $(\kappa, h)$.  
\vspace{2ex}

{\sc Theorem 4.6.} 
\begin{it}
Let  $k$ be a subfield of ${\mathbb C}$ and of $K$ which contains $\zeta_{N}$, 
and $\rho : GL_{g} \rightarrow GL_{d}$ be the rational homomorphism over $k$ 
associated with $W_{\kappa - h(1,..., 1)}$. 
Then the image of 
$$
\iota_{p} : {\cal N}_{\rho}^{\rm hol}(k) \rightarrow \overline{\cal M}_{\rho}(K)
$$
is contained in ${\cal N}_{(\kappa, h)}^{\dagger}(K)$. 
\end{it}
\vspace{2ex}

{\it Proof.} \ 
This assertion follows from Theorems 3.4 and 4.4. \ $\square$ 
\vspace{3ex}

\noindent
{\sc Department of Mathematics, 
Graduate School of Science and Engineering, 
Saga University, Saga 840-8502, Japan} 
\\
{\it E-mail:} ichikawa@ms.saga-u.ac.jp 
\vspace{2ex}

\centerline{\rule[1mm]{30mm}{0.05mm}}

\vspace{2ex}

\end{document}